\documentclass[10pt]{amsart}
\usepackage{color}
\usepackage{bbm}

\allowdisplaybreaks[4]

\definecolor{c20}{rgb}{0.,0.7,0.}
\definecolor{c30}{rgb}{0.,0.,1.}
\definecolor{c40}{rgb}{0.3,0.3,0.9}
\definecolor{c50}{rgb}{1,0,0}
\definecolor{c60}{rgb}{1,0.9,0.1}

\def\LE#1{\textcolor{c40}{#1}}

\def\LE#1{#1}

\newcommand{\tbb}[1]{{\textcolor{blue}{#1}}}
\def\tbb#1{#1}
\def\ccc#1{\textcolor{c20}{#1}}
\def\ccc#1{#1}

\def\ccP#1{\textcolor{c50}{#1}}
\def\ccP#1{#1}

\newcommand{\kb}[1]{\boldsymbol{#1}}
\newcommand{\vk}[1]{\kb{#1}}

\newcommand{\abs}[1]{\left\lvert #1 \right\rvert}

\newcommand{\E}[1]{\mathbb{E}\left(#1\right)}

\newcommand{\pk}[1]{\mathbb{P} \left( #1 \right ) }

\newcommand{\R}{\mathbb{R}}

\newcommand{\N}{\mathbb{N}}
\newcommand{\inr}{\in \R}
\newcommand{\inn}{\in \N}
\newcommand{\ldot}{,\ldots,}

\newcommand{\limit}[1]{\lim_{#1 \to   \infty}}

\newcommand{\BQN}{\begin{eqnarray}}
\newcommand{\EQN}{\end{eqnarray}}
\newcommand{\BQNY}{\begin{eqnarray*}}
\newcommand{\EQNY}{\end{eqnarray*}}

\newcommand{\BS}{\begin{sat}}
\newcommand{\ES}{\end{sat}}
\newcommand{\BT}{\begin{theo}}
\newcommand{\ET}{\end{theo}}
\newcommand{\BK}{\begin{korr}}
\newcommand{\EK}{\end{korr}}

\newcommand{\BD}{\begin{de}}
\newcommand{\ED}{\end{de}}

\newcommand{\BIT}{\begin{itemize}}
\newcommand{\EIT}{\end{itemize}}
\newcommand{\BDI}{\begin{description}}
\newcommand{\EDI}{\end{description}}

\newcommand{\BRM}{\begin{remarks}}
\newcommand{\ERM}{\end{remarks}}

\newcommand{\BEL}{\begin{lem}}
\newcommand{\EEL}{\end{lem}}

\newcommand{\BES}{\begin{sat}}
\newcommand{\EES}{\end{sat}}

\newtheorem{theo}{Theorem}[section]
\newtheorem{sat}[theo]{Proposition}
\newtheorem{de}[theo]{Definition}
\newtheorem{lem}[theo]{Lemma}

\newtheorem{korr}[theo]{Corollary}

\newtheorem{remarks}[theo]{Remarks}

\newcommand{\nelem}[1]{{Lemma \ref{#1}}}
\newcommand{\neprop}[1]{{Proposition \ref{#1}}}
\newcommand{\netheo}[1]{{Theorem \ref{#1}}}
\newcommand{\nekorr}[1]{{Corollary \ref{#1}}}

\newcommand{\COM}[1]{}

\newcommand{\QED}{\hfill $\Box$}

\def\rw{\rightarrow}

\def\IF{\infty}

\def\LT{\left}
\def\RT{\right}

\def\ooo{(1+o(1))}

\def\rw{\rightarrow}

\def\vn{\varepsilon}
\def\Var{\text{Var}}

\def\Del{\triangle}
\def\VV{\mathcal{T}}

\def\pit{\mathcal{P}}
\def\pic{\mathcal{S}}

\newcommand{\equaldis}{\stackrel{d}{=}}

\begin{document}

\title[Extremes of vector-valued Gaussian processes: exact asymptotics]
{Extremes of vector-valued Gaussian processes: exact asymptotics}

\author{Krzysztof D\c{e}bicki}
\address{Krzysztof D\c{e}bicki, Mathematical Institute, University of Wroc\l aw, pl. Grunwaldzki 2/4, 50-384 Wroc\l aw, Poland}
\email{Krzysztof.Debicki@math.uni.wroc.pl}

\author{Enkelejd  Hashorva}
\address{Enkelejd Hashorva, Department of Actuarial Science, 
University of Lausanne,\\
UNIL-Dorigny, 1015 Lausanne, Switzerland
}
\email{Enkelejd.Hashorva@unil.ch}

\author{Lanpeng Ji}
\address{Lanpeng Ji, Department of Actuarial Science, 
University of Lausanne\\
UNIL-Dorigny, 1015 Lausanne, Switzerland
}
\email{Lanpeng.Ji@unil.ch}

\author{Kamil Tabi\'{s}}
\address{Kamil Tabi\'{s}, Department of Actuarial Science, 
University of Lausanne, UNIL-Dorigny, 1015 Lausanne, Switzerland
and
Mathematical Institute, University of Wroc\l aw, pl. Grunwaldzki 2/4, 50-384 Wroc\l aw, Poland
}
\email{Kamil.Tabis@unil.ch}

\bigskip

\date{\today}
 \maketitle

{\bf Abstract:} Let $\{X_i(t),t\ge0\}, 1\le i\le n$ be mutually independent
centered  Gaussian processes with almost surely continuous sample paths.
We derive the exact asymptotics
of
\[
\pk{\exists_{t \in [0,T]}
 \forall_{i=1 \ldot n} X_i(t)> u }
\]
as $u\to\infty$, for both \LE{locally stationary} $X_i$'s and $X_i$'s with a \ccc{non-constant} generalized variance function.
Additionally, we analyze properties of multidimensional counterparts of \ccP{the
Pickands and Piterbarg} constants, that appear in the derived asymptotics.
Important by-products of \ccP{this} contribution
are the vector-process extensions of the Piterbarg inequality,
the Borell-TIS inequality, the Slepian lemma and the Pickands-Piterbarg lemma
which are the main pillars of the extremal
theory of vector-valued Gaussian processes.



\section{Introduction}

Consider a vector-valued Gaussian  process
$\{{\bf{X}}(t),t\ge 0\}$, where ${\bf{X}}(t)=(X_1(t) \ldot X_n(t))$
with $\{X_i(t),t\ge 0\}, 1\le i\le n$, $n\inn$, being \ccP{independent} centered  Gaussian processes
with almost surely (a.s.) continuous sample paths.
In this paper we focus on the asymptotic behaviour of the
probability that $\bf{X}$ enters the upper  orthant
$\{(x_1 \ldot x_n): x_i>u, i\in \{1 \ldot n\}\}$ \ccP{over a fixed time interval $[0,T]$},  i.e.,
\begin{eqnarray}
\pk{\exists_{t \in [0,T]}
 \forall_{i=1 \ldot n} X_i(t)> u }
 \label{aim}
\end{eqnarray} as $u\to\infty$.

One of important motivations
to analyze (\ref{aim}) is its
connection with the {\it conjunction} problem for Gaussian processes.
The set of conjunctions $C_{T,u}$ on \ccP{the} fixed time interval $[0,T]$
with respect to some threshold $u$ is defined as 
$$ C_{T,u}:= \{t\in [0,T]: \min_{1 \le i \le n} X_i(t) > u\} $$
see e.g., the seminal contribution \cite{MR1747100}.
One of the key properties of $C_{T,u}$, that recently focused substantial attention,
is the
probability that $C_{T,u}$ is non-empty
\BQN
 p_{T,u}:=
\pk{C_{T,u} \not=\phi}= \pk{\sup_{t\in [0,T]} \min_{1 \le i \le n} X_i(t) > u}.\label{aim2}
\EQN
\ccc{Clearly,} $p_{T,u}$ is equivalent to (\ref{aim}),
implying that one can view at (\ref{aim}) as at the probability of extremal
behaviour of {the  process}  $\{\min_{1 \le i \le n} X_i(t), t\ge 0\}$.
Typically, in applications such as the
analysis of functional magnetic resonance imaging (fMRI) data,  $X_i$'s are assumed to be real-valued Gaussian random fields.
We refer to, e.g., \cite{MR2775212, ChengXiao13,  MR1747100},
for approximations of $p_{T,u}$ in the case of smooth Gaussian random fields.
Results for non-Gaussian random fields {and general stationary processes can be found in \cite{MR2654766,DebOrderStats}}.

In the special case when $n=1$, \ccc{then} (\ref{aim})
reduces to the the tail asymptotics of supremum of a centered Gaussian process.
One of the techniques that \ccP{was}
found to be particularly successful in finding exact asymptotic  behaviour
of supremas of Gaussian processes \ccP{is} the {\it double-sum method}.
This method was originally \ccP{introduced} for the stationary case
in seminal papers of J. Pickands III \cite{PicandsB, PicandsA}.
Later, it was extended to non-stationary Gaussian processes
(and fields) {including \LE{locally stationary} Gaussian process and Gaussian process with a \ccc{non-constant} variance function.}
For a complete survey on related results we refer to \cite{Pit96,Pit20}.

The main goal of this contribution is to derive
exact asymptotics of (\ref{aim})  for  large classes of  non-stationary
Gaussian processes $X_i$'s, providing
multidimensional counterparts of the seminal
Pickands' and Piterbarg-Prishyaznyuk's results, respectively; see e.g., Theorem D2 and Theorem D3 in \cite{Pit96}.
The proofs of \ccP{our} main results are based on an extension of the double-sum
technique applied to the analysis of (\ref{aim}).
Remarkably, \ccP{the} relation between (\ref{aim}) and (\ref{aim2})
also implies the applicability of the
double-sum method
to non-Gaussian processes, \ccP{as, e.g., the process}
$\{\min_{1 \le i\le n} X_i(t), t\ge0\}$.



Interestingly, in the obtained asymptotics, there appear multidimensional
counterparts of the classical Pickands \ccP{and} Piterbarg constants
(see Sections 2 and 3). We analyze properties of these new constants in Section 3.

In the literature there are  few results on extremes of non-smooth vector-valued Gaussian processes; see  \cite{Anshin05, Debicki10, HJ14c,Yimin15} and the references therein.
In Section 5  we shall  present some extensions (tailored for our use)
of the Slepian lemma,
the Borell-TIS inequality and the Piterbarg inequality
for vector-valued Gaussian random fields.
These results are of independent interest given their crucial role in the theory
of Gaussian processes and random fields; see e.g., \cite{AdlerTaylor, DEJ13, NPV14, Pit96} and the reference\ccP{s} therein. 

{\it The organization of the paper:}
\ccP{Basic} notation and some preliminary results are presented in Section 2.
In Section 3 we analyze properties of vector-valued Pickands and Piterbarg constants.
 The main results of the paper,  concerning the asymptotics of  (\ref{aim})  for both \LE{locally stationary} $X_i$'s and $X_i$'s with a non-constant generalized variance function, are displayed in Section 4. All the proofs are relegated to Section 5.

\section{Notation and preliminaries}
We shall use some standard notation which is common when dealing with vectors.
All the operations on vectors are meant componentwise, for instance, for any given
 $ \vk{x} = (x_1,\ldots,x_n)\in \R ^n$ and $\vk{y} = (y_1,\ldots,y_n) \in \R ^n $, we write $ \vk{x} > \vk{y} $ if and only if
  $ x_i > y_i $  for all $ 1 \leq i \leq n $, {write $1/\vk{x}=(1/x_1\ccc{\ldot}1/x_n)$ if $x_i\neq 0, 1 \leq i \leq n$}, and write $\vk{x}\vk{y}=\LE{(x_1y_1 \ldot x_ny_n)}$.
Further we set  $ \vk{0}  := (0,\ldots,0)\ccP{\in\R^n} $ and $ \vk{1} : = (1,\ldots,1)\ccP{\in\R^n}$. \\
We use the notation $ f(u) = h(u)(1+o(1))$ if $ \lim_{u \to \infty} \frac{f(u)}{h(u)} = 1 $ and write $ f(u) = o(h(u)) $ if $ \lim_{u \to \infty} \frac{f(u)}{h(u)} = 0 $.
\ccP{By} $\Psi(\cdot)$ we denote  the survival function of an
$N(0,1)$ random variable, and
$\Gamma(\cdot)$ denotes the Euler Gamma function.


We shall refer to $\{\vk{X}(t), t\ge0\}$
as
a centered $n$-dimensional {\it vector-valued} Gaussian process,
where $ \vk{X}(t)=(X_1(t) \ldot X_n(t))$ with $X_i$'s being \ccP{independent} centered Gaussian processes with a.s. continuous sample paths.
Since $n$ hereafter is always fixed we shall occasionally omit
"$n$-dimensional", mentioning simply that $\vk{X}$ is a centered vector-valued Gaussian \LE{process}.
Define next
\begin{eqnarray*}
\sigma_{\vk{X}}^2(\cdot)   =   (\sigma_{X_1}^2(\cdot),\ldots,\sigma_{X_n}^2(\cdot)), \quad  R_{\vk{X}}(\cdot,\cdot)  =  (R_{X_1}(\cdot,\cdot)\ldot R_{X_n}(\cdot,\cdot)),
\end{eqnarray*}
with $\sigma_{X_i}^2(t)  =  \mathrm{Var}(X_i(t)) $ and $R_{X_i}(s,t)  =  \mathrm{Cov}(X_i(s),X_i(t))$.

Let in the following
$\{B_{i,\kappa}(t),t\in \R\}, 1\le i\le n$  be $n$ mutually independent standard fractional Brownian motions (fBm's) defined on $\R$ with \ccP{common} Hurst
index $\kappa/2\in(0,1]$, and set $\vk{B}_\kappa(t)= ( B_{1, \kappa}(t) \ldot B_{n,\kappa}(t)).$

A key step in the investigation of the tail asymptotics of supremum of Gaussian processes is the derivation of the tail asymptotic behaviour of the supremum taken over
"short intervals". For the stationary case this is achieved by the so-called Pickands lemma.
The non-stationary case is covered by the so-called Piterbarg lemma
(see \cite{DHL14Ann, DHJParisian, HJ14d} for similar terminology and related results).
Before deriving an extension of these classical results for the vector-valued \ccP{Gaussian} processes, we need to introduce some
further notation.

Let $ \{ Y(t) , t \in \R \} $ be a centered Gaussian process with a.s. continuous sample
paths such that $ Y(0) = 0 $ a.s., and let $ d \colon \R \to \R  $ be a
continuous function such that $ d(0) = 0 $. Further, denote
$S_1,S_2$ to be two non-negative constants satisfying $\max(S_1,S_2)>0$.\\
Let $ \{ X_u(t), t \in [-S_1,S_2] \}, u > 0$ be a family
of centered Gaussian processes with
a.s. continuous sample paths
that satisfies
\begin{itemize}
\item[\textbf{P1:}] $ \sigma_{X_u}^2(0) = 1 $ {for all $u$ large} and $ \lim_{u \to \infty} u^2 (1 - \sigma_{X_u}(t)) = d(t) $ uniformly with respect to $ t \in [-S_1,S_2]$;
\item[\textbf{P2:}] $ \lim_{u \to \infty} u^2 \mathrm{Var}(X_u(t)-X_u(s)) = 2 \mathrm{Var}(Y(t)-Y(s)) $ for all $ t,s \in [-S_1,S_2]$;
\item[\textbf{P3:}] there exist $ G, u_0 > 0 $ and $ \gamma \in (0,2] $ such that $ u^2 \mathrm{Var}(X_u(t)-X_u(s)) \leq G |t-s|^\gamma $ holds for all $ u \geq u_0 $ and $ s,t \in [-S_1,S_2]$.
\end{itemize}
We  write
$ X_u \in \pit(Y,d) $ if $\{X_u\}_{\ccP{u>0}}$ satisfies {\bf P1-P3}.

Introduce next some further notation which is related to vector version of the
Pickands and Piterbarg constants.
Consider $ \{ \vk{Y}(t), t \in \R \} $, with $\vk{Y}(t)=(Y_1(t) \ldot Y_n(t))$, where
$ \{ Y_i(t), t \in \R \} $ are mutually independent Gaussian processes
with a.s. continuous sample
paths such that $ Y_i(0) = 0 $ a.s., and let $ \vk{d}(t)=(d_1(t) \ldot d_n(t))$
with
$d_i(\cdot)$ being continuous  functions such that
$ d_i(0) = 0 $.
We define
\begin{eqnarray*}
\mathcal{H}_{\vk{Y} ,\vk{d} }[-S_1,S_2]   &:=&
 \int_{\R ^n} e^{\sum_{i=1}^n w_i} \mathbb{P} \Bigl( \exists_{t \in [-S_1,S_2]} \: \sqrt{2} \vk{Y} (t) - \sigma_{\vk{Y} }^2(t) - \vk{d} (t) > \vk{w} \Bigr) d \vk{w} \\
& =&
\int_{\R ^n} e^{\sum_{i=1}^n w_i} \mathbb{P} \Biggl( \sup_{t \in [-S_1,S_2]} \min_{1 \leq i \leq n} \Bigl( \sqrt{2} Y_i(t) - \sigma_{Y_i}^2(t) - d_i(t) - w_i \Bigr) > 0 \Biggr) d\vk{w}\in(0,\IF).
\end{eqnarray*}

In the special case of
$\vk{Y}(t)=\vk{B}_\kappa(t)$
being an
$n$-dimensional vector-valued fBm process with independent coordinates
we set
$$ \mathcal{H}_{\vk{B}_\kappa}(S_2): =  \mathcal{H}_{\vk{B}_\kappa,\vk{0} }[0,S_2].$$

The above defined constants play significant role in the
following multidimensional extension of the Pickands-Piterbarg lemma
(compare with, e.g., \cite{DHL14Ann, DebKo2013, Pit96}).

\BES \label{propA}
Let $ \{ {\vk{X}}_u(t), t \in [-S_1,S_2] \}, u > 0 $ be a family of centered \ccP{vector-valued} Gaussian
process  with  independent coordinates
$ X_{i,u} \in \pit(Y_i,d_i)$ for some $Y_i, d_i, 1\le i\le n$.
If $\vk{f}(\cdot)$ is an $n$-dimensional vector function such that
$ \lim_{u \to \infty} \vk{f} (u)/u = \vk{c}  > \vk{0}  $, then as $ u \to \infty $
$$ \mathbb{P} \Bigl( \exists_{t \in [-S_1,S_2]} \: {\vk{X}}_u(t) > \vk{f} (u) \Bigr)   =
\mathcal{H}_{\vk{c} \vk{Y} ,\vk{c} ^2\vk{d} }[-S_1,S_2] \prod_{i=1}^n \Psi(f_i(u)) (1+o(1)).
$$
\EES
The proof of Proposition \ref{propA} is given in Section \ref{s.proof.propA}.

Let $\vk{X}$  be a centered vector-valued Gaussian processes
with independent coordinates \ccP{ $X_i$'s which} are stationary Gaussian processes
with  unit  variance and correlation functions $r_i(\cdot), 1\le i\le n$ satisfying 
 \BQN \label{stationaryR}
r_i(t) = 1 - a_i \abs{t}^{\kappa_i} + o(\abs{t}^{\kappa_i}) \  \ t \to 0,\ \ {\rm and}\  r_i(t)< 1, \  \forall t \not=0,
\EQN
where  $\kappa_i  \in (0,2]$,
$a_i>0,$ $1\le i\le n$. Let $\kappa=\min_{1\le i\le n} \kappa_i$, and denote $\vk{a}=(a_1 1_{\{\kappa_1=\kappa\}},...,a_n 1_{\{\kappa_n=\kappa\}})$ with  $1_{\{\cdot\}}$ denoting the indicator function.
 Hereafter we write $\vk{X} \in \mathcal{S}( \vk{a}, \kappa)$ if \eqref{stationaryR} is satisfied by the vector-valued Gaussian process $\vk{X}$.

As a straightforward implication of \neprop{propA} we obtain the following corollary.



\BK \label{lemPit}
Consider a centered vector-valued stationary Gaussian process $  {\vk{X}}   \in \pic(\vk{a}, \kappa)$. 
For any $\beta_i \ge \kappa$ and
$\vk{b}(t)$ such that
$b_i(t)=
\underline{b}_i|t|^{\beta_i}1_{\{t\le0\}}
+
\overline{b}_i|t|^{\beta_i}1_{\{t>0\}}
$
$1\le i\le n$, define $ Z_i(t)= \frac{X_i(t)}{1+ b_i(t)}, t\in\R, 1\le i\le n$.
If $\vk{f}(\cdot)$ is an $n$-dimensional vector function such that
$ \lim_{u \to \infty} \vk{f} (u)/u = \vk{c} > \vk{0}  $,
then for any  non-negative constants $S_1,S_2$
satisfying  $\max(S_1,S_2)>0$ 
\BQN
 \mathbb{P} \Bigl( \exists_{t \in [-S_1u^{-2/\kappa},S_2u^{-2/\kappa}]} \: \vk{Z}(t)> \vk{f} (u) \Bigr)   =
\mathcal{H}_{\vk{c}  \sqrt{\vk{a} }\vk{B}_\kappa, \vk{c} ^2\vk{d} }[-S_1,S_2] \prod_{i=1}^n \Psi(f_i(u)) (1+o(1)) 
\EQN
holds as $u\to \IF$,
where $\vk{d}(t)=(d_1(t),...,d_n(t))$
with
$d_i(t)=b_i(t)1_{\{\beta_i=\kappa\}}$.
\EK

Next we introduce multidimensional counterparts  of \ccP{the} Pickands constant, defined as
\BQN \label{PicD}
\mathcal{H}_{\vk{C}\vk{B}_\kappa} \; := \; \lim_{S \to \infty} S^{-1} \mathcal{H}_{\vk{C}\vk{B}_\kappa}(S)
\EQN
for $ \kappa \in (0,2] $ and $ \vk{C}  \geq \vk{0}  $, $ \vk{C}  \neq \vk{0}$.
Note that
if $n=1$ {and $C_1\neq0$}, then
$\mathcal{H}_{\vk{C}\vk{B}_\kappa}=
C_1^{2/\kappa}\mathcal{H}_\kappa^*$,
where
\[\mathcal{H}_\kappa^*= \lim_{T\to\infty} T^{-1}\E {\exp\LT(\sup_{t\in[0,T]}(\sqrt{2}B_{1,\kappa}(t)-t^\kappa)\RT)}
\]
is the classical Pickands constant; see e.g., \cite{Pit96} and the recent contributions \LE{\cite{DieMik15,DikerY, Yakir}.}
{The existence and finiteness of ${\mathcal{H}_{\vk{C}\vk{B}_\kappa} }$ follow by Fekete's Lemma, since
by Lemma \ref{l.sub} displayed in Section {\ref{s.pickands},
$\mathcal{H}_{\vk{C}\vk{B}_\kappa}(S)$ is sub-additive. Furthermore, Proposition \ref{Th07} below shows that
${\mathcal{H}_{\vk{C}\vk{B}_\kappa} }>0$  }}.

Finally we introduce  multidimensional counterparts  of
Piterbarg \tbb{constants}. For $ \kappa \in (0,2]$ let
$\underline{\vk{d}}=(\underline{d}_1 \ldot \underline{d}_n)$,
$\overline{\vk{d}}=(\overline{d}_1 \ldot \overline{d}_n)$
be such that
$\sum_{i=1}^n \underline{d}_i > 0$ and
$\sum_{i=1}^n \overline{d}_i > 0$,
and let
$\vk{d}(t)=(d_1(t) \ldot d_n(t))$ with
$d_i(t)=
\underline{d}_i|t|^{\kappa}1_{\{t\le0\}}
+
\overline{d}_i|t|^{\kappa}1_{\{t>0\}}.
$
We define, for $\vk{C}\ge0$, ${\vk{C}\neq \vk{0}}$,
\BQNY
\mathcal{H}_{\vk{C}\vk{B}_\kappa}^{\underline{\vk{d}}}&:=&
 \lim_{S \to \infty}  \mathcal{H}_{\vk{C}\vk{B}_\kappa, \vk{d}}[-S,0]\nonumber\\
\mathcal{H}_{\vk{C}\vk{B}_\kappa}^{\overline{\vk{d}}}&:=&
\lim_{S \to \infty}  \mathcal{H}_{\vk{C}\vk{B}_\kappa, \vk{d}}[0,S]\nonumber\\
\mathcal{H}_{\vk{C}\vk{B}_\kappa}^{\underline{\vk{d}},\overline{\vk{d}}}&:=&
 \lim_{S \to \infty}  \mathcal{H}_{\vk{C}\vk{B}_\kappa, \vk{d}}[-S,S].
\EQNY
In Theorem \ref{ThmNS} we shall prove that the above generalized
Piterbarg constants exist and are both positive and finite.

\section{Estimates of the generalized Pickands  and  Piterbarg constants}\label{s.pickands}

In this section we provide some estimates of {the above defined}
multidimensional counterparts  of Pickands and Piterbarg constants. 
We begin with the subadditivity property of $\mathcal{H}_{\vk{C}\vk{B}_\kappa}(S)$.
\BEL\label{l.sub}
Let $ \kappa \in (0,2] $ and $ \vk{C}  \geq \vk{0}  $, $ \vk{C}  \neq \vk{0}  $.
Then for all $S \inn$
\BQN \label{subH}
 \mathcal{H}_{\vk{C} \vk{B}_\kappa}(S) \; \leq \;
 S \ \mathcal{H}_{\vk{C} \vk{B}_\kappa}(1)\in(0,\IF).
 \EQN
\EEL
The proof of Lemma \ref{l.sub} is given in Section \ref{s.proof.sub}.

Clearly, from the subadditivity of $\mathcal{H}_{\vk{C}\vk{B}_\kappa}(\cdot)$
we obtain that $\mathcal{H}_{\vk{C}\vk{B}_\kappa}$ \ccP{exists and is finite}.
In the next proposition  we confirm that $\mathcal{H}_{\vk{C}\vk{B}_\kappa}$ is strictly positive
by establishing a positive lower bound.


\BS
\label{Th07}
If $ \kappa \in (0,2] $ and $ \vk{C}  \geq \vk{0}  $, $ \vk{C}  \neq \vk{0}  $, then
\BQN\nonumber
\mathcal{H}_{\vk{C}\vk{B}_\kappa}   \geq
 \frac{\left( \sum_{i=1}^n C_i^2 \right)^{1/\kappa}}{4^{1+1/\kappa} \Gamma(1/\kappa+1)}.
\EQN
\ES

\BS
\label{Prop.lower}
For each $n\inn $ we have
\BQN\nonumber
\mathcal{H}_{\vk{B}_1} \; \leq \; n  \left( \frac{n}{n-1}  \left(2+\sqrt{\frac{2}{\pi e}}\right) \right)^{n-1} \quad
\text{ and }\quad
\mathcal{H}_{\vk{B}_2} \; \leq \; n  \left( \frac{n}{n-1} \right)^{n-1} \frac{1}{\sqrt{\pi}} ,
\EQN
where $n/(n-1)$ is set to \ccP{be} 1 for $n=1$.
\ES

We conclude this section with  lower bounds for the generalized Piterbarg constants
$\mathcal{H}_{\vk{C}\vk{B}_\kappa}^{\overline{\vk{d}}},
\mathcal{H}_{\vk{C}\vk{B}_\kappa}^{\underline{\vk{d}},\overline{\vk{d}}}$.

\BS \label{lem:low} For any $ \kappa \in (0,2]$,
$ \vk{C}  \geq \vk{0},\vk{C}  \neq \vk{0}$ and
$\underline{\vk{d}},\overline{\vk{d}}$ \LE{satisfying}
$\sum_{i=1}^n \underline{d}_i>0,\sum_{i=1}^n \overline{d}_i>0$
we have
\BQN\nonumber
\mathcal{H}_{\vk{C}\vk{B}_\kappa}^{\overline{\vk{d}}}
\ge
\ \Bigl( e\kappa \sum_{i=1}^n \tbb{\max(0,\overline{d}_i)} \Bigr)^{-1/ \kappa}
\mathcal{H}_{\vk{C}\vk{B}_\kappa } \EQN
and
\BQN\nonumber
\mathcal{H}_{\vk{C}\vk{B}_\kappa}^{\underline{\vk{d}},\overline{\vk{d}}}
\ge
\ \tbb{ 2 \left(  e\kappa\right)^{-1/ \kappa}
 \left(\sum_{i=1}^n \left(\max(0,\underline{d}_i)
 +
 \max(0,\overline{d}_i)\right)\right)^{-1/ \kappa}
\mathcal{H}_{\vk{C}\vk{B}_\kappa }.}
\EQN
\ES

We note that the lower bounds above are new even for the case $n=1$.

\section{Main Results}\label{s.main}

In this section 
we derive the asymptotics of \eqref{aim} for $\vk{X}$ with \LE{locally stationary} coordinates (see e.g., \cite{Ber74, Berman92, Hus90, Pit96}  for  \LE{locally stationary} Gaussian processes) in \netheo{Th06}
and for a large class of $\vk{X}$ with a non-constant generalized variance function in \netheo{ThmNS}. These results provide multidimensional counterparts of   Pickands theorem
 and Piterbarg-Prishyaznyuk theorem.
\subsection{Locally stationary coordinates}

{Let
 $\{\vk{X}(t),t\in[0,T]\}$  be a centered vector-valued Gaussian process
with independent coordinates $X_i$'s which are \LE{locally stationary} Gaussian processes \LE{with  a.s. continuous sample paths,    unit  variance} and correlation functions $r_i(\cdot,\cdot), 1\le i\le n$ satisfying
 \BQN \label{stationaryR0}
r_i(t,t+h) = 1 - a_i(t) \abs{h}^{\kappa_i} + o(\abs{h}^{\kappa_i}), \  \ h \to 0
\EQN
uniformly with respect to $t\in[0,T]$,
where  $\kappa_i  \in (0,2]$, and
$a_i(t), 1\le i\le n$ are positive continuous functions on $[0,T]$, 
and further
 \BQN \label{stationaryR2}
  r_i(s,t)< 1, \  \forall s,t\in[0,T]\ \mathrm{and\ } s \not=t.
\EQN
Let in  the following
$\vk{a}(t)=(a_1(t)1_{\{\kappa_1=\kappa\}},...,a_n(t)1_{\{\kappa_n=\kappa\}}),  t\in[0,T]$.
\tbb{Recall that we set $\kappa=\min_{1\le i\le n}\kappa_i$.}
}

Note that  $\vk{X}\in\pic(\vk{a}, \kappa)$ is a particular example of the above defined vector-valued Gaussian processes.


\BT \label{Th06}
Let $\vk{X}$  be a centered vector-valued Gaussian process with independent \LE{locally stationary} coordinates satisfying \eqref{stationaryR0} and \eqref{stationaryR2}.  
If $\vk{f}(\cdot)$ is an $n$-dimensional vector function such that $ \lim_{u \to \infty} \vk{f} (u)/u = \vk{c}  > \vk{0}  $,  then
\tbb{
\BQN
 \mathbb{P}  \left( \exists_{t \in [0,T]} \; {\vk{X}}(t) > \vk{f} (u) \right) \; = \;
\int_0^T \mathcal{H}_{\vk{c}  \sqrt{\vk{a}(t) }\vk{B}_\kappa} dt
 \ u^{\frac{2}{\kappa}} \prod_{i=1}^n \Psi(f_i(u)) (1+o(1)), \quad u\to \IF.
 \EQN
 }
\ET

\tbb{The special case of Theorem \ref{Th06} for $\vk{X}\in\pic(\vk{a}, \kappa)$
has been derived in \cite{DELK}}.
\LE{A straightforward comparison of Theorem \ref{Th06}}
with Theorem 1.1 in \cite{DELK} implies the following proposition.

\BS
If  $\vk{C}\ge \vk{0}$, $\vk{C}\neq \vk{0}$ and $\kappa\in(0,2]$, then
$$
\mathcal{H}_{\vk{C}\vk{B}_\kappa} =
\lim_{u \downarrow 0} \frac{1}{u}\pk{ \max_{k\ge 1}
\mathcal{Z}_{\vk{C} \vk{B}_\kappa} (uk)\le 0} \in (0,\IF)
,
$$
with
\BQN\label{mathcZ}\nonumber
\mathcal{Z}_{\vk{C} \vk{B}_\kappa} (t)
 :=  \min_{1 \le i \le n}{ \Bigl( \sqrt{2} C_i  B_{i,\kappa}(t)
- C_i^2 t^{\kappa}  + E_i\Bigr)},\ \ t\ge0, 
\EQN
where
${E_i}$'s are
mutually independent unit mean \LE{exponential random variables}  being further independent of
${B_{i, \kappa}}$'s.
\ES

\subsection{General  non-stationary coordinates}

Let $\{\vk{X}(t),t\in[0,T]\}$ be a centered vector-valued non-stationary Gaussian process with a non-constant generalized variance function.
The following set of conditions constitute\ccP{s} a vector-valued counterpart of
Piterbarg-type conditions on $X_i$'s
(see e.g., \cite{Pit96} for the original Piterbarg's
conditions imposed on   Gaussian processes or fields with a \ccc{non-constant} variance function):

{\bf Assumption I:} The following {\it generalized variance} function
$$g(t)=\sum_{i=1}^n \frac{1}{\sigma_{X_i}^2(t)}$$
 attains its minimum over $[0,T]$ at the unique point $t=t_0\in[0,T]$.

{\bf Assumption II:} There exist $\alpha_i\in(0,2]$, $a_i>0,$ $ 1\le i\le n$
such that 
$$
\mathrm{Cov}\LT(\frac{X_i(t)}{\sigma_{X_i}(t)},\frac{X_i(s)}{\sigma_{X_i}(s)}\RT)
=
1-a_i\abs{t-s}^{\alpha_i}-o(\abs{t-s}^{\alpha_i}),\ \ 1\le i\le n
$$
holds as $t,s\to t_0$.

{\bf Assumption III:} There exist some $\beta>0$,
$\underline{\vk{b}}=(\underline{b}_1 \ldot \underline{b}_n)$ and
$\overline{\vk{b}}=(\overline{b}_1 \ldot \overline{b}_n)$
such that
$$
1- \frac{\sigma_{X_i}(t_0+t)}{\sigma_{X_i}(t_0)}=
\underline{b}_{i}\abs{t}^{\beta}1_{\{t\le0\}}+
\overline{b}_{i}\abs{t}^{\beta}1_{\{t>0\}}
+o(\abs{t}^{\beta}),\ \ 1\le i\le n
$$
holds as $t\to 0$.

\def\uth{\underline{\theta}}
\def\oth{\overline{\theta}}
\LE{Note  in passing} that Assumption {III} implies that
\begin{eqnarray}
{g(t_0+t)-g(t_0)}=
2\left(\uth 1_{\{t\le0\}}
+
\oth 1_{\{t>0\}}
\right) |t|^\beta+o(\abs{t}^{\beta}),\ \ t\to0,
 \label{bound.sigma}
\end{eqnarray}
  which combined with Assumption {I} implies
$$\uth:= \sum_{i=1}^n \frac{\underline{b}_i}{\sigma_{X_i}^2(t_0)}\ge0, \quad  \oth:=\sum_{i=1}^n \frac{\overline{b}_i}{\sigma_{X_i}^2(t_0)}\ge0.$$

{\bf Assumption IV:} \LE{There exist some} positive constants $G, \gamma$ and $\rho$ such that 
%
$$
\LE{ \max_{1 \le i \le n}  \E{(X_i(t )-X_i(s ))^2}}\le G\abs{t-s}^{\gamma} 
$$
holds for  all $s,t\in(t_0-\rho,t_0+\rho)\cap[0,T]$.


\BT \label{ThmNS}
Let $\vk{X}$ be a centered vector-valued Gaussian process that satisfies Assumptions I--IV
with the parameters therein.
  Denote $\alpha=\min_{1\le i\le n}\alpha_i$, $\vk{a}=(a_11_{\{\alpha_1=\alpha\}} \ldot a_n1_{\{\alpha_n=\alpha\}})$, and let $\vk{c}=(c_1 \ldot c_n)$ with $c_i=\frac{1}{\sigma_{X_i}(t_0)}, 1\le i\le n.$
\ccc{Suppose that $\uth>0$ and $\oth>0$}. 

i) If $\alpha<\beta$, then
\BQNY
\pk{\sup_{t\in[0,T]}\min_{1\le i\le n}X_i(t)> u }=
\mathcal{H}_{\vk{c} \sqrt{\vk{a}}\vk{B}_\alpha}
\Theta\
\Gamma\LT(\frac{1}{\beta}+1\RT)
u^{\frac{2}{\alpha}-\frac{2}{\beta}}\prod_{i=1}^n\Psi\LT(c_i u \RT)\ooo, \quad u\to \IF,
\EQNY
where
$$
\Theta=
\left\{\begin{array}{lll}
 \oth ^{-\frac{1}{\beta}},& t_0=0\\
 \uth ^{-\frac{1}{\beta}}
+ \oth ^{-\frac{1}{\beta}},&  t_0\in(0,T) \\
 \uth^{-\frac{1}{\beta}},& t_0=T.
\end{array}\right.
$$
ii) If $\alpha=\beta$, then
\BQNY
\pk{\sup_{t\in[0,T]}\min_{1\le i\le n}X_i(t)> u }=
\widehat{\mathcal{H}}
\prod_{i=1}^n \Psi\LT(c_i u \RT)\ooo, \quad u\to \IF,
\EQNY
where
\BQNY
\widehat{\mathcal{H}}
=\left\{\begin{array}{lll}
\mathcal{H}_{\vk{c} \sqrt{\vk{a}}\vk{B}_\alpha}^ { \vk{c}^2\vk{\overline{b}}}   ,& t_0=0\\
\mathcal{H}_{\vk{c} \sqrt{\vk{a}}\vk{B}_\alpha}
^{\vk{c}^2\vk{\underline{b}},\vk{c}^2\vk{\overline{b}}},&  t_0\in(0,T)\\
\mathcal{H}_{\vk{c} \sqrt{\vk{a}}\vk{B}_\alpha}^ { \vk{c}^2\vk{\underline{b}}}   ,& t_0=T.
 \end{array}\right.
\EQNY
iii) If $\alpha>\beta$, \ccc{then}
\BQNY
\pk{\sup_{t\in[0,T]}\min_{1\le i\le n}X_i(t)> u }=     \prod_{i=1}^n\Psi\LT(c_i u\RT)\ooo, \quad u\to \IF.
\EQNY
\ET

{\bf Remarks}: {\it a) For $n=1$, the above theorem reduces to the classical result  for non-stationary Gaussian processes (see e.g., \cite{Pit96, MR2819240}).\\ 
b) Let $\vk{X}$ be a centered vector-valued Gaussian process \ccP{with independent coordinates $X_i$'s which are  copies of a Gaussian process $X$}, and let
  $\{X_{r:n}(t),t\ge0\}, 1\le r\le n$
be the order statistics processes of $\{X_i(t),t\ge0\}, 1\le i\le n$, i.e., we  define
 \BQNY 
X_{1:n}(t):=\max_{1 \le i \le n} X_i(t)\ge X_{2:n}(t)\ge \ldots \ge X_{n:n}(t)=\min_{1 \le i \le n} X_i(t),\ \ \ t\ge0.
\EQNY
Under the assumptions of {\netheo{Th06} or} \netheo{ThmNS}, with similar arguments as in \cite{DELK} we obtain 
$$
\pk{\sup_{t\in[0,T]} X_{r:n}(t)> u }=\frac{n!}{(n-r)! r!} \pk{\sup_{t\in[0,T]}\min_{1\le i\le r}X_i(t)> u }\ooo,\ \ u\to\IF.
$$
}

\section{Proofs}

Before proceeding to the proofs of \netheo{Th06} and  Theorem \ref{ThmNS},
we present four lemmas that will play  important \ccP{roles} in further
analysis \LE{and being} also of some independent interest.
We begin with a vector version of the Slepian lemma,
then give the vector-valued counterparts of the Borell-TIS inequality and the Piterbarg inequality,
respectively. 
Below  we write $\mathcal{T}$ for a compact set in ${\R^k, k\ge 1}$ \LE{and denote by $\abs{x}$ the Euclidean norm of $x\inr^k$} .

\BEL\label{MultiSlepian} (Slepian Lemma)
Let $ \{ {\vk{Y}}(t) , t \in \VV \} $ and $ \{ \vk{Z} (t) , t \in \VV \} $ be
two centered separable vector-valued Gaussian processes with independent coordinates.
If  for all $ s,t \in \mathcal{T} $
\begin{eqnarray*}
\sigma_{{\vk{Y}}}^2(t)  =  \sigma_{\vk{Z} }^2(t), \quad  R_{{\vk{Y}}}(t,s) \ge  R_{\vk{Z} }(t,s),
\end{eqnarray*}
then for any $ \vk{u} \in \R ^n $ we have
\BQN
\mathbb{P} \Bigl( \exists_{t \in \mathcal{T}} \: {\vk{Y}}(t) > \vk{u} \Bigr) \; \leq \; \mathbb{P} \Bigl( \exists_{t \in \mathcal{T}} \; \vk{Z} (t) > \vk{u} \Bigr).
\EQN
\EEL

\LE{{\bf Proof:}}
\COM{
Let
$T_m \subset [0,T]$ such that for each $m$ the set $T_m$ has a finite number of elements and $T_m \subset T_{m+1},m\ge 2$.
Define next
$$ X_{k,i} := X_i(t_k) , \quad  Y_{k,i} := Y_i(t_k) , \quad  \lambda_{k,i} := -u_i, \quad 1 \leq k \leq m, 1 \leq i \leq n  .$$
Note that $ \{ {\vk{X}}(t), t \in \mathcal{T}_m \} \equaldis \{ - {\vk{X}}(t)  t \in \mathcal{T}_m \} $ implying that
\BQNY
\mathbb{P} \Bigl( \exists_{t \in \mathcal{T}_m} \: {\vk{X}}(t) > \vk{u} \Bigr) & = & \mathbb{P} \Bigl( \exists_{t \in \mathcal{T}_m} \: - {\vk{X}}(t) > \vk{u} \Bigr)\\
& = &\mathbb{P} \Biggl( \bigcup_{k=1}^m \bigcap_{i=1}^n \Bigl\{ - X_i(t_k) > u_i \Bigr\} \! \Biggr) \\[1ex]
& = & 1 - \mathbb{P} \Biggl( \bigcap_{k=1}^m \bigcup_{i=1}^n \Bigl\{ - X_i(t_k) \leq u_i \Bigr\} \! \Biggr) \\
&= & 1 - \mathbb{P} \Biggl( \bigcap_{k=1}^m \bigcup_{i=1}^n \Bigl\{ X_{k,i} \geq \lambda_{k,i} \Bigr\} \! \Biggr) \: .
\EQNY
By the mutual independence of coordinates we have for all $ 1 \leq k,l \leq m, 1 \leq i \neq j \leq n $
\[ \E{X_i(t_k) X_j(t_l)} \; = \; \E{X_{k,i} X_{l,j}} \; = \;
\E{Y_{k,i} Y_{l,j}}  = \E{Y_i(t_k) Y_j(t_l)} .
\]
Furthermore, in view of our assumptions for all $ 1 \leq k \neq l \leq m, 1 \leq i \leq n $
\[ \E{X_i^2(t_k)}  = \E{X_{k,i}^2 }= \E{ Y_{k,i}^2 }= \E{Y_i^2(t_k)} , \]
\[ \E{ X_i(t_k) X_i(t_l) } = \E{X_{k,i} X_{l,i}}\geq \E{Y_{k,i} Y_{l,i}} = \E{ Y_i(t_k) Y_i(t_l)}.
\]
Hence, the assumptions of Gordon's inequality (see \cite{Gordon}) are satisfied, consequently
\begin{equation} \label{lab106} \mathbb{P} \biggl( \sup_{t \in \mathcal{T}_m} \min_{1 \leq i \leq n} \Bigl( X_i(t) - u_i \Bigr) > 0 \biggr) \; = \; \mathbb{P} \Bigl( \exists_{t \in \mathcal{T}_m} \: {\vk{X}}(t) > \vk{u} \Bigr) \; \leq \; \mathbb{P} \Bigl( \exists_{t \in \mathcal{T}_m} \: \vk{Y} (t) > \vk{u} \Bigr) \end{equation}
for any finite $ \mathcal{T}_m $.
}
The claim for 
 any finite set $\mathcal{T}$ follows by a direct application of Gordon's inequality (see \cite{Gordon}). If
 $\mathcal{T}$ is a given compact set {of $\R^k$},  then the proof can be easily established  using standard arguments that make use of the separability assumption; \ccP{see e.g., \cite{AdlerTaylor}.} \QED
 \COM{
Let
$\mathcal{T }$ be a finite set with $m$ elements, and we write $\mathcal{T} =\{t_1,\cdots,t_m\}$. 
Further define
$$ Y_{i,j} := Y_j(t_i) , \quad  Z_{i,j} := Z_j(t_i), \quad 1 \leq i \leq m, 1 \leq j \leq n.$$
Note that $ \{ {\vk{Y}}(t), t \in \mathcal{T} \} \equaldis \{ - {\vk{Y}}(t),  t \in \mathcal{T}\} $. Thus we have 
\BQNY
\mathbb{P} \Bigl( \exists_{t \in \mathcal{T}} \: {\vk{Y}}(t) > \vk{u} \Bigr) & = & \mathbb{P} \Bigl( \exists_{t \in \mathcal{T}} \: - {\vk{Y}}(t) > \vk{u} \Bigr)\\
& = &\mathbb{P} \Biggl( \bigcup_{i=1}^m \bigcap_{j=1}^n \Bigl\{ - Y_j(t_i) > u_j \Bigr\} \! \Biggr) \\
& = &\mathbb{P} \Biggl( \bigcup_{i=1}^m \bigcap_{j=1}^n \Bigl\{ Y_{i,j} \le - u_j \Bigr\} \! \Biggr).
\EQNY
Next, we verify the  assumptions of Gordon's inequality. First we have
$$
\E{Y_{i,j}^2 }= \E{(Y_j(t_i))^2}= \sigma^2_{Y_j}(t_i)= \sigma^2_{Z_j}(t_i) = \E{(Z_j(t_i))^2} = \E{ Z_{i,j}^2 }
 $$
holds for $1 \leq i \leq m, 1 \leq j \leq n$. Further, by the mutual independence of coordinates we have, for any $1 \leq i \leq m, 1 \leq j,k \leq n,$
$$
\E{ Y_{i,j} Y_{i,k} } = \E{Y_j(t_i) Y_k(t_i)}= \E{ Z_j(t_i) Z_k(t_i)}=\E{Z_{i,j} Z_{i,k}}.
$$
Moreover for any $ i \neq l, 1 \leq i, l \leq m, 1 \leq j,k  \leq n$  we have,
\BQNY
\E{Y_{i,j} Y_{l,k}}= \E{Y_j(t_i) Y_k(t_l)}\ge \E{ Z_j(t_i) Z_k(t_l)}=\E{Z_{i,j} Z_{l,k}}.
 \EQNY
Therefore, the assumptions of Gordon's inequality  are satisfied. Consequently
\begin{equation} \label{lab106} \mathbb{P} \biggl( \sup_{t \in \mathcal{T}} \min_{1 \leq i \leq n} \Bigl( Y_i(t) - u_i \Bigr) > 0 \biggr) \; = \;
\mathbb{P} \Bigl( \exists_{t \in \mathcal{T}} \: {\vk{Y}}(t) > \vk{u} \Bigr) \; \leq \; \mathbb{P} \Bigl( \exists_{t \in \mathcal{T}} \: \vk{Z} (t) > \vk{u} \Bigr)
\end{equation}
holds for any finite set $ \mathcal{T}$.

We now consider the general compact $\mathcal{T}$. For each $m > 0$ let $\mathcal{T}_m$ be a finite subset of $\mathcal{T}$ such that
$\mathcal{T}_m \subset \mathcal{T}_{m+1}$ and $\mathcal{T}_m$ increases to a dense subset of $\mathcal{T}$. By separability,
we have $\cup_{m\ge 1} \mathcal{T}_m= \mathcal{T}$ and thus 
\[ \sup_{t \in \mathcal{T}_m} \min_{1 \leq i \leq n} \Bigl( Y_i(t) - u_i \Bigr) \; \stackrel{\textrm{a.s.}}{\longrightarrow} \; \sup_{t \in \mathcal{T}} \min_{1 \leq i \leq n} \Bigl( Y_i(t) - u_i \Bigr) \]
and
\[ \sup_{t \in \mathcal{T}_m} \min_{1 \leq i \leq n} \Big( Z_i(t) - u_i \Bigr) \; \stackrel{\textrm{a.s.}}{\longrightarrow} \; \sup_{t \in \mathcal{T}} \min_{1 \leq i \leq n} \Bigl( Z_i(t) - u_i \Bigr) \:  \]
hold as $m\to \IF$. Moreover, since the last convergence are monotone, the claims follows by 
 (\ref{lab106}).  \QED

}

Set in the following ${\tau^2_\VV}=\inf_{t\in\VV}\sum_{i=1}^n \frac{1}{\sigma_{X_i}^2(t)}$.

\BEL\label{GBorell} (Borell-TIS inequality)
Let \ccc{$\{\vk{X}(t),  t\in\VV\}$} be a centered vector-valued Gaussian process
with independent coordinates which have a.s. continuous sample paths.  If $\tau_\VV>0$,  then
there exists some positive constant $\mu$ such that for $u>\mu$
\BQNY
\pk{\exists_{t\in\VV}\vk{X}(t)>u\vk{1}}\le \exp\LT(-\frac{(u-\mu)^2}{2}{\tau^2_\VV}\RT).
\EQNY
\EEL
{\bf Proof:}
It follows that
\BQN\label{eq:BorellY}
\pk{\exists_{t\in\VV}\vk{X}(t)>u\vk{1}}&\le&\pk{\sup_{t\in\VV} Y(t)>u },
\EQN
where (set $A(t)=\sum_{i=1}^n\prod_{j=1,j\neq i}^n\sigma_{Z_j}^2(t), t\in \VV$)
$$
Y(t)=\sum_{i=1}^n \Biggl( \frac{\prod_{j=1,j\neq i}^n\sigma_{Z_j}^2(t)}{A(t)}\Biggr)X_i(t),\ \ \ t\in\VV.
$$
Since further
\BQNY
\Var\LT(Y(t)\RT)=\LT(\sum_{i=1}^n \frac{1}{\sigma_{X_i}^2(t)}\RT)^{-1}
\EQNY
the claim follows from the Borell-TIS inequality for one-dimensional Gaussian processes (e.g., \cite{AdlerTaylor}) with
$$\mu=\E{\sup_{t\in \VV} Y(t)  }<\IF$$
and thus the proof is complete.  \QED

\BEL\label{GPiter} (Piterbarg inequality)
\ccc{Under the conditions of \nelem{GBorell}, if further Assumption IV holds},
then  for all $u$ large
\BQN\label{IPIT}
\pk{\exists_{t\in\VV}\vk{X}(t)>u\vk{1}}
\le C \text{mes}(\VV) \ u^{\frac{2}{\nu}-1} \exp\LT(-\frac{u^2}{2}{\tau^2_\VV}\RT),
\EQN
where $C$ is some positive constant not depending on $u.$
\EEL

{\bf Proof:} We use  the same notation as in the proof of \nelem{GBorell}. In the light of \eqref{eq:BorellY} and Theorem 8.1 in \cite{Pit96}, it suffices to show that
\BQN \label{eq:regu}
\E{(Y(t) -Y(s) )^2}\le L \abs{t-s}^\nu,\ \ \forall  s,t \in \VV
\EQN
holds for some positive \LE{constant $L$,} \ccc{which is a direct consequence of Assumption IV}. \QED

\COM{

{\bf Remarks:}
 {\it 
a)  The set $\VV$ in the last lemmas  can be  more general.
For instance, \nelem{MultiSlepian} holds if we assume
$(\mathcal{T}, d)$ to be  a separable metric space,
and  \nelem{GBorell} holds if we assume $\{\vk{X}(t),  t\in\VV\}$ to be a
centered vector-valued Gaussian process  with independent coordinates which have a.s. bounded sample paths on a general parameter space  $\VV$. \\
\COM{Since the proof of Piterbarg inequality strongly relies on the Borell-TIS  inequality it is possible also to allow $\mathcal{T}$ to be a general set. For example
 $\mathcal{T}$ can be a compact subset of $\R^k,k \ge 1$.\\
b)  The claim in \eqref{IPIT} is valid if one of $\nu_1$ or $\nu_2$
 is equal to infinity. This means that the claim in \eqref{IPIT} is valid if we assume only say the condition \eqref{eq:ZG1} (or similarly only condition \eqref{eq:ZG2}). \\
}
b) Assume that
$\mathcal{G}=\LT\{t\in \VV: t=\text{arg}\inf \sum_{i=1}^n \frac{1}{\sigma_{Z_i}^2(t)}\RT\}$
is a finite set. Define
$\mathcal{G}_\vn= \underset{t\in \mathcal{G}} \bigcup ([t-\vn,t+\vn]\cap\VV)$
for any small positive $\vn$. In view of the proof of  Theorem 8.1 in \cite{Pit96},
the claim of \nelem{GPiter} still holds if 
\eqref{eq:ZG2} is valid for all $t,s\in \mathcal{G}_\vn$ \ccP{with} some small positive $\vn$.
}\\
}

The  last lemma  below concerns the asymptotics of a probability of double events; it  is crucial when dealing with the double sum term in the proof of our main results.

\BEL \label{Lem20} Consider a  centered
vector-valued stationary Gaussian process $  {\vk{X}}   \in \pic(\vk{a}, \kappa)$.
 Suppose that for those $X_i$'s with $\kappa_i=\kappa$ there exists some global constant $\vn>0$ such that
\[
 1 - \frac{a_i}{2} t^{\kappa} \geq R_{X_i}(t,0) \geq 1 - 2a_i t^{\kappa}  
 \]
holds for all $ t \in [0,\vn]$.
Assume further that $\vk{f}(\cdot), \vk{h}(\cdot)$ are two  $n$-dimensional vector functions such that $\lim_{u \to \infty} \vk{f} (u)/u = \vk{c} _1> \vk{0}$ and $ \lim_{u \to \infty} \vk{h}(u)/u = \vk{c} _2> \vk{0}$. 
Then there exist two positive constants $ F, G $ such that
 for all $ t_0 > S > 1 $
\BQNY \lefteqn{\mathbb{P}  \left( \exists_{t \in [0,S]u^{-2/\kappa}} \; {\vk{X}}(t) > \vk{f} (u) , \exists_{t \in [t_0,t_0+S]u^{-2/\kappa}} \; {\vk{X}}(t) > \vk{h}(u) \right)}\\
& \leq & F S^{2n} \exp \left( - G (t_0-S)^\kappa \right) \prod_{i=1}^n \Psi\LT(\frac{f_i(u)+h_i(u)}{2}\RT)
\EQNY
holds for all $u$ large. 
\EEL
{\bf Proof:}
First note that if $\kappa_i=\kappa$, in view of the proof of Lemma 6.3 in
\cite{Pit96} (or   Lemma 5 in \cite{Michna09}) we obtain that
\BQNY
&&\pk{ \sup_{t \in [0,S]u^{-2/\kappa}} X_i(t) > f_i(u) ,
\sup_{t \in [t_0,t_0+S]u^{-2/\kappa}} X_i(t) > h_i(u)}   \\
 && \leq   F_i S^2 \exp \left( - G_i (t_0-S)^\kappa \right) \Psi\LT(\frac{f_i(u)+h_i(u)}{2}\RT)
\EQNY
holds with some positive constants $F_i, G_i$.
Further, if $\kappa_i>\kappa$, then there \LE{exist} some positive constant $L$ and sufficiently small $\vn_1>0$ such that
$$
r_{i}(t)\ge e^{-L t^\kappa}
$$
is valid for all $t\in[0,\vn_1]$.\\
 Let $\{\xi(t),t\ge0\}$ be a stationary
Gaussian process with a.s. continuous sample paths and correlation
function $r_{\xi}(t)=e^{-L t^\kappa}, t\ge0$. By the Slepain lemma
(cf. Theorem C.1 in \cite{Pit96} or \nelem{MultiSlepian}) we have 
\BQNY
 &&\mathbb{P}  \left( \sup_{t \in [0,Su^{-2/\kappa}]} X_i(t) > f_i(u) ,
 \sup_{t \in [t_0,t_0+S]u^{-2/\kappa}} X_i(t) > h_i(u) \right)\\
&& \leq \mathbb{P}  \left( \sup_{t \in [0,Su^{-2/\kappa}]} X_i(t) > \max(f_i(u),h_i(u)) \right) \\
 &&\leq\pk{\sup_{t \in [0,Su^{-2/\kappa}]} \xi(t) > \frac{f_i(u)+h_i(u)}{2}}. 
 \EQNY
Consequently, \LE{the} Pickands lemma (cf. Lemma D.1 in \cite{Pit96} or \nekorr{lemPit}) implies
 \BQNY
 &&\mathbb{P}  \left( \sup_{t \in [0,Su^{-2/\kappa}]} X_i(t) > f_i(u) ,
 \sup_{t \in [t_0,t_0+S]u^{-2/\kappa}} X_i(t) > h_i(u) \right)\\
&& \leq  F_i S^2 \Psi\LT(\frac{f_i(u)+h_i(u)}{2}\RT)\le  F_i S^2 \exp \left( - G_i (t_0-S)^\kappa \right)\Psi\LT(\frac{f_i(u)+h_i(u)}{2}\RT)
 \EQNY
for all $u$ sufficiently large, with $G_i = 0$ and some $F_i > 1$. Moreover since in view of the independence of $ X_i$'s 
\BQNY
&&\mathbb{P}  \left( \exists_{t \in [0,S]u^{-2/\kappa}} \; {\vk{X}}(t) > \vk{f} (u) , \exists_{t \in [t_0,t_0+S]u^{-2/\kappa}} \; {\vk{X}}(t) > \vk{h}(u) \right)\\
&& \leq \; \mathbb{P}  \left( \bigcap_{i=1}^n \left\{ \sup_{t \in [0,S]u^{-2/\kappa}} X_i(t) > f_i(u) \right\} , \bigcap_{i=1}^n \left\{ \sup_{t \in [t_0,t_0+S]u^{-2/\kappa}} X_i(t) > h_i(u) \right\} \right) \\
&& = \; \prod_{i=1}^n \mathbb{P}  \left( \sup_{t \in [0,S]u^{-2/\kappa}} X_i(t) > f_i(u) , \sup_{t \in [t_0,t_0+S]u^{-2/\kappa}} X_i(t) > h_i(u) \right)
\EQNY
the claim follows by choosing $ F = \prod_{i=1}^n F_i > 0 $, $ G = \sum_{i=1}^n G_i > 0 $. This completes the proof. \QED



\subsection{Proof of Proposition \ref{propA}}\label{s.proof.propA}

The idea of the proof is based on a multidimensional modification of the proof of
Theorem D.1 in \cite{Pit96}. We shall present only the main steps that lead to the claim.
For all $ u > 0 $ we have
\begin{eqnarray*}
\lefteqn{
 \mathbb{P} \left( \exists_{t \in [0,T]} \; {\vk{X}}_u(t) > \vk{f}(u) \right)
=\int_{\mathbbm{R}^n} \prod_{i=1}^n  \left( \frac{1}{\sqrt{2 \pi}} e^{- v_i^2 / 2} \right)
\mathbb{P}  \left( \exists_{t \in [0,T]} \; {\vk{X}}_u(t) > \vk{f}(u) \;
\Big{|} \; {\vk{X}}_u(0) = \vk{v}  \right) d\vk{v} \nonumber}\\
&=& \prod_{i=1}^n  \left( \Psi(f_i(u)) \right)
\int_{\mathbbm{R}^n} e^{\sum_{i=1}^n  \left( w_i - w_i^2 / (2f_i^2(u)) \right)}
 \mathbb{P}  \left( \exists_{t \in [0,T]} \; {\vk{X}}_u(t) > \vk{f}(u) \; \Big{|} \; {\vk{X}}_u(0) = \vk{f}(u) - \frac{\vk{w} }{\vk{f}(u)} \right) d\vk{w} .
\label{lab51}
\end{eqnarray*}

Consider  the family
$ {\vk{\chi}}_u(t) = \left(\chi_{1,u}(t),\ldots,\chi_{n,u}(t)\right) $ indexed by $u$, where
$$\chi_{i,u}(t)
:=
\left( f_i(u)  \left( X_{i,u}(t) - f_i(u) \right) + w_i \;
\Big{|} \; X_{i,u}(0) = f_i(u) - \frac{w_i}{f_i(u)} \right),
$$
and observe that
\[ \mathbb{P}  \left( \exists_{t \in [0,T]} \; {\vk{X}}_u(t) > \vk{f}(u) \; \Big{|} \; {\vk{X}}_u(0) = \vk{f}(u) - \frac{\vk{w} }{\vk{f}(u)} \right)
= \mathbb{P}  \left( \exists_{t \in [0,T]} \; \vk{\chi}_u(t) > \vk{w}  \right). \]

By {\bf P1-P2} for any $w\inr$
\[
\E {\chi_{i,u}(t)}\to -
\mathrm{Var}(Y_i(t))-d_i(t), \ \ u\to\infty
\]
holds uniformly with respect to $t\in[0,T]$.
Moreover,
\[
\E{(\chi_{i,u}(t)-\chi_{i,u}(s))^2}\to
2\mathrm{Var}(Y_i(t)-Y_i(s)), \ \ u\to\infty
\]
\ccP{for all} $s,t\in[0,T]$.
Hence
\[ \lim_{u\to\infty}
 \mathbb{P}  \left( \exists_{t \in [0,T]} \; \vk{\chi}_u(t) > \vk{w}  \right)
=
 \mathbb{P}  \left( \exists_{t \in [0,T]} \; \left(\sqrt{2}\vk{Y}(t) -
\mathrm{Var}(\vk{Y}(t))-\vk{d}(t)\right) > \vk{w}  \right)
\]
for each $w\in \R$. The remaining part of the proof follows line-by-line the same reasoning as
the corresponding proof of Lemma D.1 in  \cite{Pit96}, {where {\bf P3} \tbb{is} used for the tightness of $\chi_{i,u}$'s;
\LE{see also Proposition 9.7 in \cite{Pit20} and Lemma 2 in \cite{dieker2005extremes}}.} This completes the proof.
\QED

\subsection{Proof of Lemma \ref{l.sub}}\label{s.proof.sub}

It suffices to suppose that in Corolarry \ref{lemPit}  we have
{$b_1(t)=\cdots=b_n(t)=0$} (so $Z_i(\cdot)=X_i(\cdot)$ is stationary)
and note that
\begin{eqnarray*}
\mathbb{P}  \left( \exists_{t \in [0,Su^{-2/\kappa}]} \; {\vk{X}}(t) > \vk{f} (u) \right) & \leq & \sum_{k=1}^S \mathbb{P}  \left( \exists_{t \in [k-1,k]u^{-2/\kappa}} \; {\vk{X}}(t) > \vk{f}(u) \right) \\[1ex]
& = & S \mathbb{P}  \left( \exists_{t \in [0,u^{-2/\kappa}]} \; {\vk{X}}(t) > \vk{f}(u) \right)
\end{eqnarray*}
is valid for all $u>0$.
\QED

\subsection{Proof of Proposition \ref{Th07}}\label{p.Th07}

The idea of the proof is based on a multidimensional modification of a technique
developed in Lemma 16 and  Corollary 17 in \cite{DRolski}
and in Lemma 7 in \cite{Shao}.
For a fixed $ a > 0 $ and a positive integer  $ N $, using Bonferroni's inequality,
we obtain
\BQNY
 \mathcal{H}_{\vk{C}\vk{B}_\kappa }(aN) &= & \int_{\mathbbm{R}^n} e^{\sum_{i=1}^n w_i} \mathbb{P}  \left( \exists_{t \in [0,aN]} \left\{ \sqrt{2} \vk{C}\vk{B}_\kappa(t) - \vk{C} ^2 t^\kappa > \vk{w} \right\} \right) d\vk{w} \\
& \geq & \int_{\mathbbm{R}^n} e^{\sum_{i=1}^n w_i} \mathbb{P}  \left( \exists_{1 \leq k \leq N} \left\{ \sqrt{2} \vk{C}\vk{B}_\kappa(ak) - \vk{C} ^2 (ak)^\kappa > \vk{w} \right\} \right) d\vk{w} \\
& \geq &\sum_{k=1}^N \int_{\mathbbm{R}^n} e^{\sum_{i=1}^n w_i} \mathbb{P}  \left( \sqrt{2} \vk{C}\vk{B}_\kappa(ak) - \vk{C} ^2 (ak)^\kappa) > \vk{w} \right) d\vk{w} \\
&& - \sum_{k=1}^{N-1} \sum_{l=k+1}^N  \int_{\mathbbm{R}^n}  e^{\sum_{i=1}^n w_i} \mathbb{P}  \left( \sqrt{2} \vk{C}  (\vk{B}_\kappa(ak) + \vk{B}_\kappa(al)) - \vk{C} ^2 ((ak)^\kappa + (al)^\kappa) > 2 \vk{w} \right) d\vk{w} \\
& =& \sum_{k=1}^N \int_{\mathbbm{R}^n} e^{\sum_{i=1}^n w_i} \: \prod_{i=1}^n \mathbb{P}  \left( \sqrt{2} C_iB_{i,\kappa}(ak) - C_i^2 (ak)^\kappa > w_i \right) d\vk{w} \label{lab29}\\
&& - \sum_{k=1}^{N-1} \sum_{l=k+1}^N  \int_{\mathbbm{R}^n}  e^{\sum_{i=1}^n w_i} \LE{\prod_{i=1}^n  \mathbb{P}}\left( \sqrt{2} C_i(B_{i,\kappa}(ak) + B_{i,\kappa}(al)) -  C_i^2 ((ak)^\kappa + (al)^\kappa) > 2 w_i \right) d\vk{w} \\
& = & \sum_{k=1}^N \prod_{i=1}^n \int_{\mathbbm{R}} e^{w_i} \mathbb{P}  \left( \sqrt{2} C_iB_{i,\kappa} (ak) - C_i^2 (ak)^\kappa > w_i \right) dw_i \\
&& - \sum_{k=1}^{N-1} \sum_{l=k+1}^N \prod_{i=1}^n \int_{\mathbbm{R}} e^{w_i} \mathbb{P}  \left( \frac{ \sqrt{2} C_i(B_{i,\kappa} (ak) +  B_{i,\kappa} (al)) - C_i^2 ((ak)^\kappa + (al)^\kappa) }{2} > w_i \right) dw_i \\
&= & \sum_{k=1}^N \prod_{i=1}^n \mathbb{E} \Biggl(\exp  \left( \sqrt{2} C_i B_{i,\kappa} (ak) - C_i^2 (ak)^\kappa \right)\Biggr) \\
&& \label{lab150}  - \sum_{k=1}^{N-1} \sum_{l=k+1}^N \prod_{i=1}^n \mathbb{E} \Biggl(\exp  \left( \frac{ \sqrt{2} C_i (B_{i,\kappa} (ak) +  B_{i,\kappa} (al)) - C_i^2 ((ak)^\kappa + (al)^\kappa) }{2} \right) \Biggr).
\EQNY 
 Since $ \exp  \left( \sqrt{2} C_i B_{i,\kappa} (ak) - C_i^2 (ak)^\kappa \right) $ is log-normal distributed 
\begin{eqnarray*}
\mathbb{E} \Biggl( \exp  \left( \sqrt{2} C_i B_{i,\kappa} (ak) - C_i^2 (ak)^\kappa \right)\Biggr)   &=&   1 , \\
\E{\exp  \left( \frac{ \sqrt{2} C_i (B_{i,\kappa} (ak) +  B_{i,\kappa} (al)) - C_i^2 ((ak)^\kappa + (al)^\kappa) }{2} \right) }  &=&   \exp  \left( - C_i^2(a(l-k))^\kappa / 4 \right)
\end{eqnarray*}
implying   
\begin{eqnarray*}
\mathcal{H}_{\vk{C}\vk{B}_\kappa}(aN) & \geq & N \; - \; \sum_{k=1}^{N-1} \sum_{l=k+1}^N \exp  \left( - \frac{\sum_{i=1}^n C_i^2(a(l-k))^\kappa}{4} \right) \\[1ex]
& \geq & N  \left( 1 - \sum_{k=1}^{N} \exp  \left( - \frac{C^2(ak)^\kappa}{4} \right) \right) ,
\end{eqnarray*}
where $ C^2 = \sum_{i=1}^n C_i^2 $. From the definition of $ \mathcal{H}_{\vk{C}\vk{B}_\kappa} $, for any $ a > 0 $
\BQNY
\nonumber \mathcal{H}_{\vk{C}\vk{B}_\kappa} & = & \lim_{N \to \infty} \frac{\mathcal{H}_{\vk{C}\vk{B}_\kappa}(aN)}{aN} \; \geq \; \frac{1}{a}  \left( 1 - \sum_{k=1}^{\infty} \exp  \left( - \frac{C^2 a^\kappa}{4} k^\kappa \right) \right) \\
\nonumber & \geq & \frac{1}{a}  \left( 1 - \int_{0}^{\infty} \exp  \left( - \frac{C^2 a^\kappa}{4} x^\kappa \right) dx \right) \\
 & = & \frac{1}{a}  \left( 1 - \frac{\Gamma(1/\kappa)}{\kappa  \left( C^2 a^\kappa / 4 \right)^{1/\kappa}} \right)\\
 & = & \frac{1}{a}  \left( 1 - \frac{1}{a} \frac{\Gamma(1/\kappa)}{\kappa  \left( C^2 / 4 \right)^{1/\kappa}} \right) .
\EQNY
The maximum over $ a > 0 $ of $ f(a) = \frac{1}{a}  \left( 1 - \frac{c}{a} \right) $ is attained at $ a^* = 2c $ with $ f(a^*) = \frac{1}{4c} $. Consequently, setting $ c = \frac{\Gamma(1/\kappa)}{\kappa  \left( C^2 / 4 \right)^{1/\kappa}} $ we obtain
\[ \mathcal{H}_{\vk{C}\vk{B}_\kappa} \; \geq \; \frac{\kappa  \left( C^2 / 4 \right)^{1/\kappa}}{4 \Gamma(1/\kappa)} \;  \]
establishing the  claim.
\QED
\subsection{Proof of Proposition \ref{Prop.lower}}
In view of Lemma 2.2 and Lemma 2.3 in \cite{DeKisow}, we have for any $ T > 0 $
\begin{eqnarray*}
\mathcal{H}_{B_1}(T) & = &  \left( 2 + T \right) \Phi  \left( \sqrt{T/2} \right) + \sqrt{T/\pi} \exp  \left( - T/4 \right) \; \leq \; 2 + \sqrt{\frac{2}{\pi e }} + T,  \\[1ex]
\mathcal{H}_{B_2}(T) & = & 1 + \frac{T}{\sqrt{\pi}}.
\end{eqnarray*}
Hence the case $n=1$ is clear. Next, for $n\ge 2$, from 
\ccP{the subadditivity} of
 $ \mathcal{H}_{\mathbf{B}_\kappa}(\cdot) $ and  the independence of $ B_{i,\kappa}(\cdot) $    we have
\[ \mathcal{H}_{\mathbf{B}_\kappa} \; \ccP{=} \; \inf_{T > 0 } \frac{\mathcal{H}_{\mathbf{B}_\kappa}(T)}{T} \; \leq \; \inf_{T > 0 } \frac{(\prod_{i=1}^n \mathcal{H}_{B_{i,\kappa}}(T))}{T} \; = \; \inf_{T > 0 } \frac{(\mathcal{H}_{B_\kappa}(T))^n}{T} \; . \]
Therefore, for $ \kappa = 1 $, $ \mathcal{H}_{\mathbf{B}_\kappa} \leq \min_{x > 0} \frac{(c+x)^n}{x} $ with $ c = 2 + \sqrt{\frac{2}{\pi e }}$, and the minimum is attained at $ x^* = \frac{c}{n-1} $.
For $ \kappa = 2 $, $ \mathcal{H}_{\mathbf{B}_\kappa} \leq \min_{x > 0} \frac{(1+cx)^n}{x} $ with $ c = \frac{1}{\sqrt{\pi}} $, and the minimum is attained for $ x^* = \frac{1}{(n-1)c} $.
 This completes the proof.
\QED

\subsection{Proof of Proposition \ref{lem:low}}
\label{p.lem:low}
It is sufficient to \ccc{show} the proof for
$\mathcal{H}_{\vk{C}\vk{B}_\kappa}^{\overline{\vk{d}}}$.
By definition for any $T>0$ we have
\BQNY
\lim_{S \to \infty}  \mathcal{H}_{\vk{C}\vk{B}_\kappa, \vk{d}}[0,S] &\ge &
\mathcal{H}_{\vk{C}\vk{B}_\kappa, \vk{d}}[0,T] \\
&=&
 \int_{\R ^n} e^{\sum_{i=1}^n w_i} \mathbb{P}
 \Bigl( \exists_{t \in [0,T]} \: \sqrt{2} \vk{C}\vk{B}_\kappa(t)
   - \vk{C}^2\sigma_{\vk{B}_\kappa}^2(t) - \overline{\vk{d}}  t^{\kappa} > \vk{w} \Bigr) d \vk{w} \\
& \ge&
 \int_{\R ^n} e^{\sum_{i=1}^n w_i} \mathbb{P} \Bigl( \exists_{t \in [0,T]} \: \sqrt{2} \vk{C}\vk{B}_\kappa(t) - \vk{C}^2\sigma_{\vk{B}_\kappa}^2(t)  > \vk{w}+
  \tbb{\max(\vk{0},\overline{\vk{d}})}T^{\kappa} \Bigr) d \vk{w} \\
 &=&  e^{-\ccc{d^+} T^{\kappa} }\mathcal{H}_{\vk{C}\vk{B}_\kappa, \vk{0}}[0,T], \quad d^+:=\sum_{i=1}^n\max(0,\overline{d}_i).
\EQNY
\tbb{Since}
$\mathcal{H}_{\vk{C}\vk{B}_\kappa, \vk{0}}[0,T]=\mathcal{H}_{\vk{C}\vk{B}_\kappa}(T)$ is subadditive,
Fekete's Lemma implies
\BQNY
\lim_{S \to \infty}  \mathcal{H}_{\vk{C}\vk{B}_\kappa, \vk{d}}[0,S]
&\ge & \sup_{T> 0} \Biggl(e^{- \ccc{d^+} T^{\kappa}} \mathcal{H}_{\vk{C}\vk{B}_\kappa, \vk{0}}[0,T] \Biggr)\\
&\ge &   \sup_{T>0} Te^{ -\ccc{d^+} T^{\kappa}}  \inf_{T> 0} \frac{\mathcal{H}_{\vk{C}\vk{B}_\kappa}(T)}{T}\\
&= &   \sup_{T>0}\Bigl( Te^{ -\ccc{d^+}T^{\kappa}} \Bigr)  \limit{T}  \frac{\mathcal{H}_{\vk{C}\vk{B}_\kappa }(T) }{T}\\
&= &  \LT( \ccc{d^+} e \kappa \RT)^{-1/ \kappa}  \mathcal{H}_{\vk{C}\vk{B}_\kappa }
\EQNY
establishing the proof. \QED

\subsection{Proof of Theorem \ref{Th06}}
{The complete proof consists of two steps. In Step 1 we show the claim for ${\vk{X}}$ with stationary coordinates, and then in Step 2 we show the proof for  ${\vk{X}}$ with \LE{locally stationary} coordinates.}

{\bf Step 1. Stationary coordinates, i.e., ${\vk{X}}\in \pic(\vk{a}, \kappa)$. }

First let $ S > 1 $ and denote for $u>0$
$$ \Delta_{k} = [kSu^{-\frac{2}{\kappa}},(k+1)Su^{-\frac{2}{\kappa}}], k\inn_0,\ \ \ \  N(u) = \left\lfloor T u^{2/\kappa}S^{-1} \right \rfloor+1. $$
Here $\lfloor \cdot\rfloor$ denotes the ceiling function.
By Bonferroni's inequality and the stationarity of $ {\vk{X}}$ for sufficiently large $ u $ we have
\BQNY
\pk{\exists_{t \in [0,T]} \; {\vk{X}}(t) > \vk{f}(u)} & \leq & \sum_{k=0}^{N(u)} \pk{ \exists_{t \in \Delta_{k}} \; {\vk{X}}(t) > \vk{f}( u) } \\ & = & N(u) \pk{ \exists_{t \in \Delta_0} \; {\vk{X}}(t) > \vk{f}(u) }.
\EQNY
Thus, by \nekorr{lemPit}  we obtain
\BQN\label{eq:Supper}
\limsup_{u \to \infty} \frac{ \mathbb{P}  \left( \exists_{t \in [0,T]} \; {\vk{X}}(t) > \vk{f} (u) \right) }{ T u^{2/\kappa} \prod_{i=1}^n \Psi(f_i(u)) } \; \leq \; \frac{\mathcal{H}_{\vk{c}  \sqrt{\vk{a} }\vk{B}_\kappa}(S)}{S}.
\EQN
Again by Bonferroni's inequality
\BQNY
  \mathbb{P}  \left( \exists_{t \in [0,T]} \; {\vk{X}}(t) > \vk{f} (u) \right)
  \geq   \sum_{k=0}^{N(u) - 1} \mathbb{P}  \left( \exists_{t \in \Delta_{k}} \; {\vk{X}}(t) > \vk{f} (u) \right) -\Sigma(u)
  \EQNY
holds, where
$$
\Sigma(u)= \sum_{0 \leq k < l\le N(u)} \mathbb{P}  \left( \exists_{t \in \Delta_{k}} \; {\vk{X}}(t) > \vk{f} (u) , \exists_{t \in \Delta_{l}} \; {\vk{X}}(t) > \vk{f} (u) \right).
 $$
Similarly to  the proof of \eqref{eq:Supper} we obtain
\BQN\label{eq:Slower}
\liminf_{u \to \infty} \frac{ \sum_{k=0}^{N(u) - 1} \mathbb{P}  \left( \exists_{t \in \Delta_{k}} \; {\vk{X}}(t) > \vk{f} (u) \right) }{ T u^{2/\kappa} \prod_{i=1}^n \Psi(f_i(u)) } \; \ge \; \frac{\mathcal{H}_{\vk{c}  \sqrt{\vk{a} }\vk{B}_\kappa}(S)}{S}.
\EQN
Next we shall focus on the double sum term $\Sigma(u)$. We choose some small positive $\vn$ such that the assumptions in \nelem{Lem20} are satisfied. We divide $\Sigma(u)$ into three parts, say, $\Sigma_1(u)$  the sum over indexes $k,l$ such that $(l-k-1)Su^{2/\kappa}>\vn$,  $\Sigma_2(u)$ the sum over indexes $l > k + 1$ and $(l-k-1)Su^{2/\kappa}\le\vn$, and $\Sigma_3(u)$ the sum over indexes $l = k+ 1$.

For the summand of $\Sigma_1(u)$ similarly as in the proof of Lemma \ref{Lem20} we have
\BQNY
 \mathbb{P}  \left( \exists_{t \in \Delta_k} {\vk{X}}(t) > \vk{f} (u) , \exists_{t \in \Delta_l} {\vk{X}}(t) > \vk{f} (u) \right) &\le&
 \prod_{i=1}^n \mathbb{P}  \left( \sup_{t \in \Delta_k } X_i(t) > f_i(u) , \sup_{t \in \Delta_l} X_i(t) > f_i(u) \right)\\
 & \leq &\prod_{i=1}^n \mathbb{P}  \left( \sup_{(t,s) \in \Delta_k\times\Delta_l}  X_i(t) + X_i(s) > 2f_i(u) \right).
 \EQNY
Further since $r_i(t)<1$ for any $t\not=0$, then 
$$
\delta_i=\max_{(t,s) \in \Delta_k\times\Delta_l} r(s-t)<1,
$$
which yields
$$
\Var(X_i(t) + X_i(s))=2(1+r_i(s-t))<2(1+\delta_i)< 4.
$$
Therefore from the Borell-TIS inequality   for all sufficiently large $u$
\BQNY
 \mathbb{P}  \left( \exists_{t \in \Delta_k} {\vk{X}}(t) > \vk{f} (u) , \exists_{t \in \Delta_l} {\vk{X}}(t) > \vk{f} (u) \right) &\le&
 \prod_{i=1}^n  \exp\left(-\frac{(f_i(u)-m_i)^2}{1+\delta_i}  \right)
 \EQNY
holds for some positive constants  $m_i, 1\le i\le n$. Consequently,
\BQN\label{eq:SSum1}
\limsup_{u \to \infty} \frac{ \Sigma_1(u)  }{ T u^{2/\kappa} \prod_{i=1}^n \Psi(f_i(u)) } =0.
\EQN

For the summand of $\Sigma_2(u)$, we get from \nelem{Lem20} that for all sufficiently large $u$
\BQNY
 \mathbb{P}  \left( \exists_{t \in \Delta_k} {\vk{X}}(t) > \vk{f} (u) , \exists_{t \in \Delta_l} {\vk{X}}(t) > \vk{f} (u) \right)  \le
  F S^{2n} \exp \left( - G ((l-k-1)S)^\kappa \right) \prod_{i=1}^n \Psi\LT( f_i(u) \RT)
 \EQNY
holds with some positive constants $F, G$. Thus, for sufficiently large $u$
\BQNY
 \Sigma_2(u)  \leq  (N(u)+1) \sum_{l=1}^{\infty} F  S^{2n} \exp \left( - G   (lS)^\kappa \right) \prod_{i=1}^n \Psi(f_i(u)) 
\EQNY
is valid. Note that for any  $ \theta, G > 0 $ and $ S > (\theta G/2)^{-1/\theta} $ we have
\[
\sum_{k=1}^{\infty} e^{-G(kS)^\theta} \; \leq \; 2e^{-GS^\theta} .
\]
 Consequently, {for large enough $S$}
\begin{equation} \label{eq:SSum2}
\limsup_{u \to \infty} \frac{\Sigma_2(u)}{T u^{2/\kappa}\prod_{i=1}^n \Psi(f_i(u))} \; \leq \; 2 F S^{2n-1} e^{-GS^\kappa}.
\end{equation}

For the summand of $\Sigma_3(u)$, by the stationarity of $\vk{X}$ (set ${\vk{X}}_{u}(t) = {\vk{X}} (tu^{-2/\kappa})$)
we have
\BQNY
&&  \mathbb{P}  \left( \exists_{t \in \Delta_k} \; {\vk{X}}(t) > \vk{f} (u) , \exists_{t \in \Delta_l} \; {\vk{X}}(t) > \vk{f} (u) \right)  \\
&&= \mathbb{P}  \left( \exists_{t \in [0,S]} \; {\vk{X}}_u(t) > \vk{f} (u) , \exists_{t \in [S,2S]} \; {\vk{X}}_u(t) > \vk{f} (u) \right)\\
&&=\mathbb{P}  \left( \exists_{t \in [0,S]} \; {\vk{X}}_u(t) > \vk{f} (u) , \left\{ \exists_{t \in [S,S+\sqrt{S}]} \; {\vk{X}}_u(t) > \vk{f} (u)\right\}\cup \left\{\exists_{t \in [S+\sqrt{S},2S]} \; {\vk{X}}_u(t) > \vk{f} (u) \right\} \right) \\
&&\le \mathbb{P}  \left( \exists_{t \in [0,S]} \; {\vk{X}}_u(t) > \vk{f} (u) , \exists_{t \in [S+\sqrt{S},2S+\sqrt{S}]} \; {\vk{X}}_u(t) > \vk{f} (u) \right) + \mathbb{P}  \left( \exists_{t \in [S,S+\sqrt{S}]} \; {\vk{X}}_u(t) > \vk{f} (u) \right).
\EQNY
Applying Lemma \ref{Lem20} with $ t_0 = S + \sqrt{S} $ and Pickands lemma (\ccP{see
\nekorr{lemPit}}) to the last two terms above, respectively, we obtain that for sufficiently large $u, S$
\BQNY
\lefteqn{\mathbb{P}  \left( \exists_{t \in \Delta_k} \; {\vk{X}}(t) > \vk{f} (u) , \exists_{t \in \Delta_l} \; {\vk{X}}(t) > \vk{f} (u) \right)  }\\
& \leq & F S^{2n} \exp \left( - G \sqrt{S^\kappa} \right) \prod_{i=1}^n \Psi(f_i(u)) + F \mathcal{H}_{\vk{c}  \sqrt{\vk{a} }\vk{B}_\kappa}(\sqrt{S}) \prod_{i=1}^n \Psi(f_i(u))  \\
 & \overset{\eqref{subH}}\le & \label{lab98}  F S^{2n} \exp \left( - G \sqrt{S^\kappa} \right) \prod_{i=1}^n \Psi(f_i(u)) +  F \mathcal{H}_{\vk{c}  \sqrt{\vk{a} }\vk{B}_\kappa}(1) \sqrt{S} \prod_{i=1}^n \Psi(f_i(u))   \\
&\leq & F_1  \left( S^{2n} \exp \left( - G \sqrt{S^\kappa} \right) + \sqrt{S} \right) \prod_{i=1}^n \Psi(f_i(u)),
\EQNY
with some constant $ F_1 > 0 $. 
Therefore
\[ \Sigma_3(u) \; \leq \; (N(u)+1) F_1  \left( S^{2n} \exp \left( - G \sqrt{S^\kappa} \right) + \sqrt{S} \right) \prod_{i=1}^n \Psi(f_i(u)) \]
for sufficiently large $ u $, and thus
\begin{equation} \label{eq:SSum3}
\limsup_{u \to \infty} \frac{\Sigma_3(u)}{T u^{2/\kappa} \prod_{i=1}^n \Psi(f_i(u))} \; \leq \; F_1  \left( S^{2n-1} \exp \left( - G \sqrt{S^\kappa} \right) +  S^{-\frac{1}{2}} \right) .
\end{equation}
 Consequently, it follows from (\ref{eq:Supper}--\ref{eq:SSum3}) that for any sufficiently large $S_1,S_2$
\BQN\label{eq:Hca}
 \frac{\mathcal{H}_{\vk{c}  \sqrt{\vk{a} }\vk{B}_\kappa}(S_1)}{S_1} &\ge &\limsup_{u \to \infty} \frac{ \mathbb{P}  \left( \exists_{t \in [0,T]} \; {\vk{X}}(t) > \vk{f} (u) \right) }{ T u^{2/\kappa} \prod_{i=1}^n \Psi(f_i(u)) }\nonumber \\
& \ge& \liminf_{u \to \infty} \frac{ \mathbb{P}  \left( \exists_{t \in [0,T]} \; {\vk{X}}(t) > \vk{f} (u) \right) }{ T u^{2/\kappa} \prod_{i=1}^n \Psi(f_i(u)) }\label{eq:HHa}\\
& \ge& \frac{ \mathcal{H}_{\vk{c}  \sqrt{\vk{a} }\vk{B}_\kappa}(S_2) }{S_2}  -   F_1  \left( S_2^{2n-1} \exp \left( - G  \sqrt{S_2^\kappa} \right) + S_2^{-\frac{1}{2}} \right) - 2 F S_2^{2n-1} e^{-GS_2^\kappa}.\nonumber
\EQN
Hence, the claim of the theorem follows  from \eqref{eq:HHa}  by letting $S_1,S_2\to\IF$.

{For Step 2,  we \tbb{point} out that  a close observation of \eqref{eq:Hca} shows that
\BQN\label{eq:Hu}
\frac{\mathcal{H}_{\vk{c}  \sqrt{\vk{a} }\vk{B}_\kappa}(S_1)}{S_1}\to \mathcal{H}_{\vk{c}  \sqrt{\vk{a} }\vk{B}_\kappa},\ \ {S_1\to\IF}
\EQN
uniformly with respect to $\vk{a}\in[\vk{\nu},\vk{\mu}]:=\prod_{i=1}^n[\nu_i,\mu_i]$, with $\vk{\nu}\ge\vk{0}$, $ \vk{\nu}\neq \vk{0}$ and $\nu_i<\mu_i<\IF, 1\le i\le n.$}

{\bf Step 2. Locally stationary coordinates.}

We  consider only the case where $\kappa=\kappa_1=\cdots=\kappa_n$; the same approach applies for the general case.
It follows from \eqref{stationaryR0} that for any $\vn>0$ there is some small $\delta_0>0$ such that for all $1\le i\le n$
\BQN \label{eq:ait}
 (1-\vn)a_i(t) \abs{h}^{\kappa} \le 1 - r_i(t,t+h) \le (1+\vn)a_i(t) \abs{h}^{\kappa}
\EQN
hold for all $t,t+h \in[0,T]$ satisfying $\abs{h}\le \delta_0$. Now let $\lambda\in(0,\delta_0)$ be any small constant and denote $\lambda_k=k\lambda, k\inn_0$. Clearly
\BQNY
\sum_{k=0}^{\lfloor T /\lambda \rfloor+1} \pk{\exists_{t \in [\lambda_k, \lambda_{k+1}]} \; {\vk{X}}(t) > \vk{f} (u)}&\ge& \pk{\exists_{t \in [0,T]} \; {\vk{X}}(t) > \vk{f} (u)}\\
&\ge& \sum_{k=0}^{\lfloor T /\lambda \rfloor} \pk{\exists_{t \in [\lambda_k, \lambda_{k+1}]} \; {\vk{X}}(t) > \vk{f} (u)}-\Sigma_4(u),
\EQNY
with
$$
\Sigma_4(u)=\sum_{0\le k<l\le \lfloor T /\lambda \rfloor}\pk{\exists_{t \in [\lambda_k, \lambda_{k+1}]} \; {\vk{X}}(t) > \vk{f} (u), \exists_{t \in [\lambda_l, \lambda_{l+1}]} \; {\vk{X}}(t) > \vk{f} (u)}.
$$
Next, for any fixed $k\inn_0$, define centered stationary Gaussian processes $\{\xi_i^{\vn_\pm}(t), t\ge0\}$ with unit variance and correlation functions
\BQNY
r_{\xi_i^{\vn_\pm}}(t)=\exp(- (1\pm\vn)a_i(\lambda_k) \abs{t}^{\kappa}), t\ge0,\ \ \ 1\le i\le n,
\EQNY
and let $\vk{\xi}^{\vn_\pm}(t)=(\xi_1^{\vn_\pm}(t),\cdots,\xi_n^{\vn_\pm}(t)), t\ge0$.  In view of \eqref{eq:ait} and \nelem{MultiSlepian} we have
\BQNY
  \pk{\exists_{t \in [\lambda_k, \lambda_{k+1}]} \; {\vk{\xi}^{\vn_-}}(t) > \vk{f} (u)}\le   \pk{\exists_{t \in [\lambda_k, \lambda_{k+1}]} \; {\vk{X}}(t) > \vk{f} (u)}
   \le   \pk{\exists_{t \in [\lambda_k, \lambda_{k+1}]} \;\vk{\xi}^{\vn_+}(t) > \vk{f} (u)}.
\EQNY
Then applying the results in Step 1 for vector-valued stationary Gaussian process we conclude that for $\lambda$ sufficiently small
\BQN \label{eq:lsu}
\limsup_{u\to\IF}\frac{\sum_{k=0}^{\lfloor T /\lambda \rfloor+1} \pk{\exists_{t \in [\lambda_k, \lambda_{k+1}]} \; {\vk{X}}(t) > \vk{f} (u)}}{u^{\frac{2}{\kappa}} \prod_{i=1}^n \Psi(f_i(u)) }&\le&
 \sum_{k=0}^{\lfloor T /\lambda \rfloor+1}   \mathcal{H}_{\vk{c}  \sqrt{(1+\vn)\vk{a}(\lambda_k) }\vk{B}_\kappa}\  \lambda\nonumber\\
 &\le & (1+\vn)\int_{0}^T \mathcal{H}_{\vk{c}  \sqrt{(1+\vn)\vk{a}(t) }\vk{B}_\kappa}dt,
\EQN
where the last inequality follows from the fact that $\mathcal{H}_{\vk{c} \sqrt{(1+\vn)\vk{a}(t) }\vk{B}_\kappa}$ is continuous with respect to $t\in[0,T]$ which is due to \eqref{eq:Hu} and some elementary derivations. Similarly, we have for $\lambda$ sufficiently small
\BQN \label{eq:lsl}
\liminf_{u\to\IF}\frac{\sum_{k=0}^{\lfloor T /\lambda \rfloor} \pk{\exists_{t \in [\lambda_k, \lambda_{k+1}]} \; {\vk{X}}(t) > \vk{f} (u)}}{u^{\frac{2}{\kappa}} \prod_{i=1}^n \Psi(f_i(u)) }
\ge (1-\vn)\int_{0}^T \mathcal{H}_{\vk{c}  \sqrt{(1-\vn)\vk{a}(t) }\vk{B}_\kappa}dt.
\EQN
Furthermore, similar arguments as in the proof of Theorem 7.1 in \cite{Pit96} show that
\BQNY
\limsup_{u\to\IF}\frac{\Sigma_4(u)}{u^{\frac{2}{\kappa}} \prod_{i=1}^n \Psi(f_i(u)) }=0.
\EQNY
Consequently, the claim follows by letting $\vn\to0$ in \eqref{eq:lsu} and \eqref{eq:lsl}.
This completes the proof.
\QED


\subsection{Proof of \netheo{ThmNS}}
We only give the proof for the case that $t_0\in(0,T)$, {$\alpha=\alpha_1=\cdots=\alpha_n$ and $\overline{\vk{b}}>\vk{0},  \underline{\vk{b}}>\vk{0}$.}
The proofs of the other cases follow by similar arguments and are therefore omitted.

Let $\delta(u)=(\ln u/u)^{2/\beta}$, and denote {$D_u=[ -\delta(u),   \delta(u)]$} for $u$ large. In the following,   all   formulas are meant for large enough $u$. With these notation we have
\BQNY
 P_{1}(u):=\pk{\sup_{t\in (t_0+D_u)}\min_{1\le i\le n}X_i(t)> u}&\le &  \pk{\sup_{t\in[0,T]}\min_{1\le i\le n}X_i(t)> u } \\
& \le &P_{1}(u)+\pk{\sup_{t\in[0,T]/(t_0+D_u)}\min_{1\le i\le n}X_i(t)> u }\\
& =: &P_{1}(u)+P_{2}(u).
\EQNY
Next, we shall derive the exact asymptotics of $P_{1}(u)$ as $u\rw\IF$, and show that 
\BQN\label{eq:P2}
P_{2}(u)=o(P_{1}(u)),\ \ \ \ \ u\rw\IF
\EQN
implying thus
 \BQNY
\pk{\sup_{t\in[0,T]}\min_{1\le i\le n}X_i(t)> u }=P_{1}(u)(1+o(1)),\ \ \ \ \ u\rw\IF.
\EQNY
Now we focus on the the asymptotics of $P_{1}(u)$ as $u\rw\IF$. For any small enough $\vn>0$ define
\BQNY
Z_i^{\vn_{\pm}}(t)=\frac{\sigma_{X_i}(t_0)}{1+d_i^{\vn_{\mp}}(t)}
\eta_i^{\vn_{\pm}}(t),\ \  {t\in \R,}\ \ \ 1\le i\le n,
\EQNY
where
\BQN\label{eq:d}
d_i^{\vn_{\mp}}(t)=(\underline{b}_i\mp \vn )|t|^\beta 1_{\{t\le 0\}}
+
(\overline{b}_i\mp \vn )|t|^\beta 1_{\{t> 0\}},\ \ t\in\R,
\EQN
and $\{\eta_i^{\vn_{\pm}}(t), t\in\R\}$ are centered stationary Gaussian processes with unit variance and correlation functions
$$
r_{\eta_i^{\vn_{\pm}}}(t)=\exp\LT(-a_i^{\vn_\pm} \abs{t} ^\alpha\RT),\   \ccP{t\ge0},\ \ \ a_i^{\vn_\pm}= a_i\pm \vn.
$$
In view of Assuptions II--III and \nelem{MultiSlepian}, we have that for any small enough $\vn>0$
\BQN \label{eq:P1}
\pk{\sup_{t\in D_u}\min_{1\le i\le n}Z_i^{\vn_-}(t)> u}\le P_1(u)\le  \pk{\sup_{t\in D_u}\min_{1\le i\le n}Z_i^{\vn_+}(t)> u}
\EQN
holds for all $u$ sufficiently large. In the following,
we shall show that the above upper and lower bounds for $P_1(u)$
are asymptotically equivalent as $u\to\IF$ and $\vn\to 0$.

Next we introduce some notation.
Let $T_1$ be any positive constant.  For the case that $\alpha\le \beta$, we can split the interval $D_u$ into several sub-intervals of side lengths $T_1u^{-2/\alpha}$. Specifically, let
$$
\Del_{k}= \LT[ kT_1u^{-\frac{2}{\alpha}},  (k+1)T_1u^{-\frac{2}{\alpha}}\RT],\ \ \ k\in \mathbb{Z},\ \ N(u)=\LT\lfloor T_1^{-1} (\ln u)^{\frac{2}{\beta}}u^{\frac{2}{\alpha}-\frac{2}{\beta}} \RT\rfloor+1
$$
 and note that
\BQNY
 \bigcup _{k=-N(u)+1}^{N(u)-1}  \Del_{k } \subset D_u \subset  \bigcup _{k=-N(u)}^{N(u)}   \Del_{k }.
\EQNY
We deal with the three cases
i) $\alpha<\beta$,  ii) $\alpha=\beta$ and  iii) $\alpha>\beta$ one-by-one, using different techniques.

\underline{Case i) $\alpha<\beta$}:  {\it Upper bound.} Using  Bonferroni inequality we have
\BQNY
P_1(u)\le  \pk{\sup_{t\in D_u}\min_{1\le i\le n}Z_i^{\vn_+}(t)> u}\le \sum_{k=-N(u)}^{N(u)}\pk{\sup_{t\in \Del_k}\min_{1\le i\le n}Z_i^{\vn_+}(t)> u}
\EQNY
and
\BQNY
\sum_{k=0}^{N(u)}\pk{\sup_{t\in \Del_k}\min_{1\le i\le n}Z_i^{\vn_+}(t)> u}&\le&
\sum_{k=0}^{N(u)}\pk{\exists_{t\in \Del_k}\vk{\eta}^{\vn_+}(t)> \vk{f}^{\vn_-}(k,u)}\\
&=&
\sum_{k=0}^{N(u)}\pk{ \exists_{t\in [0,T_1u^{-\frac{2}{\alpha}}]} \vk{\eta}^{\vn_+}(t)> \vk{f}^{\vn_-}(k,u)},
\EQNY
where $\vk{\eta}^{\vn_+}(t)=(\eta_1^{\vn_{\pm}}(t),\cdots,\eta_n^{\vn_{\pm}}(t)), t\ge0,$  $\vk{f}^{\vn_\pm}(k,u)=(f_1^{\vn_\pm}(k,u),\cdots,f_n^{\vn_\pm}(k,u))$ with
$$
f_i^{\vn_\pm}(k,u)=
\frac{1}{\sigma_{X_i}(t_0)}
\left(1+(\underline{b}_i\pm\vn)(|k|T_1u^{-\frac{2}{\alpha}})^\beta1_{\{k\le0\}}+
(\overline{b}_i\pm\vn)(kT_1u^{-\frac{2}{\alpha}})^\beta1_{\{k>0\}}\right)u,\ \ 1\le i\le n. 
$$
Recall that we set $\vk{c}=(c_1,\cdots,c_n)$ with  $c_i=\frac{1}{\sigma_{X_i}(t_0)}, 1\le i\le n.$
Applying \nekorr{lemPit}  we obtain
\BQNY
\sum_{k=0}^{N(u)}\pk{ \exists_{t\in [0,T_1u^{-\frac{2}{\alpha}}]} \vk{\eta}^{\vn_+}(t)> \vk{f}^{\vn_-}(k,u)}= \mathcal{H}_{\vk{c} \sqrt{\vk{a}^{\vn_+}}\vk{B}_\alpha}(T_1)\sum_{k=0}^{N(u)}\prod_{i=1}^n\Psi(f_i^{\vn_-}(k,u))\ooo
\EQNY
as $u\to\IF$.
Since further, with
$\oth_{\vn\pm}:= \sum_{i=1}^n c_i^2( \overline{b}_i\pm\vn)$,
\BQNY
\sum_{k=0}^{N(u)}\prod_{i=1}^n\Psi(f_i^{\vn_-}(k,u))&=&
\sum_{k=0}^{N(u)}\prod_{i=1}^n\LT(\frac{1}{\sqrt{2\pi}f_i^{\vn_-}(k,u)}\exp\LT(-\frac{(f_i^{\vn_-}(k,u))^2}{2}\RT)\RT)\ooo\\
&=&(2\pi)^{-\frac{n}{2}}\LT(\prod_{i=1}^n\sigma_{X_i}(t_0)\RT) u^{-n}  \sum_{k=0}^{N(u)}\exp\LT(-\sum_{i=1}^n\frac{( 1+ 2(\overline{b}_i-\vn)(kT_1u^{-\frac{2}{\alpha}})^\beta)  u^2}{2 \sigma_{X_i}^2(t_0)}\RT)\ooo\\
&=&(2\pi)^{-\frac{n}{2}}\LT(\prod_{i=1}^n\sigma_{X_i}(t_0)\RT) u^{-n} \exp\LT(-\frac{u^2}{2}g(t_0)\RT) \sum_{k=0}^{N(u)}\exp\LT(-
\oth_{\vn-} (kT_1u^{\frac{2}{\beta}-\frac{2}{\alpha}})^\beta\RT)\ooo\\
&=&T_1^{-1}(2\pi)^{-\frac{n}{2}}\LT(\prod_{i=1}^n\sigma_{X_i}(t_0)\RT) u^{\frac{2}{\alpha}-\frac{2}{\beta}-n}  \exp\LT(-\frac{u^2}{2}g(t_0)\RT) \int_{0}^\IF\exp\LT(-\oth_{\vn-} x^\beta\RT)dx\ooo\\
\EQNY
we conclude that
\BQNY
\sum_{k=0}^{N(u)}\pk{\sup_{t\in \Del_k}\min_{1\le i\le n}Z_i^{\vn_+}(t)> u}&\le& \frac{\mathcal{H}_{\vk{c} \sqrt{\vk{a}^{\vn_+}}\vk{B}_\alpha}(T_1)}{T_1 \oth_{\vn-}^{1/\beta}} 
\ccc{\varphi^*}(u)\ooo,
\EQNY
where
\BQNY
\ccc{\varphi^*}(u)&:=&(2\pi)^{-\frac{n}{2}}\LT(\prod_{i=1}^n\sigma_{X_i}(t_0)\RT)
\Gamma\LT(\frac{1}{\beta}+1\RT)
u^{\frac{2}{\alpha}-\frac{2}{\beta}-n}
 \exp\LT(-\frac{u^2}{2}g(t_0)\RT)\\
&=&\Gamma\LT(\frac{1}{\beta}+1\RT)u^{\frac{2}{\alpha}-\frac{2}{\beta}}\prod_{i=1}^n\Psi(c_iu)\ooo.
\EQNY
Similarly, we can find  \ccc{an} upper bound for
$\sum_{k=-N(u)}^{0}\pk{\sup_{t\in \Del_k}\min_{1\le i\le n}Z_i^{\vn+}(t)> u}$,
which leads to
\BQN
\sum_{k=-N(u)}^{N(u)}\pk{\sup_{t\in \Del_k}\min_{1\le i\le n}Z_i^{\vn_+}(t)> u}
&\le&
 \frac{\mathcal{H}_{\vk{c} \sqrt{\vk{a}^{\vn_+}}\vk{B}_\alpha}(T_1)}{T_1}
\Theta_{\vn-}
 \ccc{\varphi^*}(u)
\ooo,
\label{eq:upp2}
\EQN
where
$
\Theta_{\vn-}=
\left( \oth_{\vn-} \right)^{-1/\beta}+
 \left(\uth_{\vn-}  \right)^{-1/\beta},
$ with
$\uth_{\vn\pm}= \sum_{i=1}^n c_i^2 ( \underline{b}_i\pm \vn)$.

{\it Lower bound.} Applying again the Bonferroni inequality we have
\BQNY
P_1(u)\ge  \pk{\sup_{t\in D_u}\min_{1\le i\le n}Z_i^{\vn_-}(t)> u}\ge \sum_{k=-N(u)+1}^{N(u)-1}\pk{\sup_{t\in \Del_k}\min_{1\le i\le n}Z_i^{\vn_-}(t)> u}-\Sigma_1(u),
\EQNY
where
$$
\Sigma_1(u)=\underset{-N(u)\le k<l\le N(u)}{\sum\sum }\pk{\sup_{t\in \Del_k}\min_{1\le i\le n}Z_i^{\vn_-}(t)> u, \sup_{t\in \Del_l}\min_{1\le i\le n}Z_i^{\vn_-}(t)> u }.
$$
With similar arguments as for the derivation of \eqref{eq:upp2} we obtain that
\BQN
\sum_{k=-N(u)+1}^{N(u)-1}\pk{\sup_{t\in \Del_k}\min_{1\le i\le n}Z_i^{\vn_-}(t)> u}
&\ge&
\frac{\mathcal{H}_{\vk{c} \sqrt{\vk{a}^{\vn_-}}\vk{B}_\alpha}(T_1)}{T_1}
\Theta_{\vn+}
\ccc{\varphi^*}(u)
\ooo,\label{eq:low2}
\EQN
where
$\Theta_{\vn+}=
\left(\oth_{\vn+} \right)^{-1/\beta}+
 \left(\uth_{\vn+}  \right)^{-1/\beta}
$.

Next we consider the double sum term $\Sigma_1(u)=:\Sigma_2(u)+\Sigma_3(u)$ where $\Sigma_2(u)$   is the sum over indexes $l = k+ 1$, and $\Sigma_3(u)$ is the sum over indexes $l > k + 1$.  It follows that
\BQNY
\Sigma_2(u)&=&\sum_{ k=-N(u)}^{ N(u)} \pk{\sup_{t\in \Del_k}\min_{1\le i\le n}Z_i^{\vn_-}(t)> u}+\sum_{ k=-N(u)}^{ N(u)}\pk{ \sup_{t\in \Del_{k+1}}\min_{1\le i\le n}Z_i^{\vn_-}(t)> u }\\
&&-\sum_{ k=-N(u)}^{ N(u)} \pk{\sup_{t\in \Del_k\cup\Del_{k+1}}\min_{1\le i\le n}Z_i^{\vn_-}(t)> u}.
\EQNY
Thus, we have from \eqref{eq:upp2} and \eqref{eq:low2} that
\BQN \label{eq:sig2}
\limsup_{u\to\IF}\frac{\Sigma_2(u)}{\ccc{\varphi^*}(u)}&\le&
2 \Theta_{\vn-} \frac{\mathcal{H}_{\vk{c} \sqrt{\vk{a}^{\vn_+}}\vk{B}_\alpha}(T_1)}{T_1}
-2 \Theta_{\vn+} \frac{
\mathcal{H}_{\vk{c} \sqrt{\vk{a}^{\vn_-}}\vk{B}_\alpha}(2T_1) }{2T_1}.
\EQN
Moreover
\BQNY
\Sigma_3(u) 
&=&\sum_{ k=0}^{ N(u)}\sum_{l= k+2}^{N(u)} \pk{\sup_{t\in \Del_k}\min_{1\le i\le n}Z_i^{\vn_-}(t)> u, \sup_{t\in \Del_{l}}\min_{1\le i\le n }Z_i^{\vn_-}(t)> u }\\
&& +\sum_{ k=-N(u)}^{-1}\sum_{l= k+2}^{N(u)} \pk{\sup_{t\in \Del_k}\min_{1\le i\le n}Z_i^{\vn_-}(t)> u, \sup_{t\in \Del_{l}}\min_{1\le i\le n }Z_i^{\vn_-}(t)> u }\\
&=:&\Sigma_{3,1}(u)+\Sigma_{3,2}(u)
\EQNY
An application of \nelem{Lem20} gives that
\BQNY
\Sigma_{3,1}(u)
&\le&\sum_{ k=0}^{ N(u)}\sum_{l= k+2}^{N(u)} \pk{\exists_{t\in \Del_k}\vk{\eta}^{\vn_-}(t)> \vk{f}^{\vn_+}(k,u), \exists_{t\in \Del_l}\vk{\eta}^{\vn_-}(t)> \vk{f}^{\vn_+}(l,u) }\\
&\le& \sum_{ k=0}^{ N(u)}\sum_{l= k+2}^{N(u)} FT_1^{2n} \exp\LT(-G ((l-k-1)T_1)^\alpha\RT)  \prod_{i=1}^n \Psi\LT(\frac{f_i^{\vn_+}(k,u)+ f_i^{\vn_+}(l,u)}{2}\RT)
 \EQNY
holds with some positive constants $F, G$ for any $u$ sufficiently large. Using the same
reasoning as in \eqref{eq:upp2} and noting that $\oth_{\vn\pm}>0, \uth_{\vn\pm}>0$ for sufficiently small $\vn$, we conclude that
\BQNY
 \limsup_{u\to\IF}\frac{\Sigma_{3,1}(u)}{\ccc{\varphi^*}(u)}
\le  F_1 T_1^{2n-1}\sum_{l= 1}^\IF \exp\LT(-G (lT_1)^\alpha\RT),
\EQNY
with some $F_1>0.$
Similar arguments apply also for $\Sigma_{3,2}(u)$ and thus we have
\BQN\label{eq:sig3}
\limsup_{T_1\to\IF}\limsup_{u\to\IF}\frac{\Sigma_3(u)}{\ccc{\varphi^*}(u)}
=0.
\EQN
Consequently, by letting $\vn\to0$ and $T_1\to\IF$ we obtain from (\ref{eq:upp2}--\ref{eq:sig3}) that
 \BQNY
P_1(u) = \mathcal{H}_{\vk{c} \sqrt{\vk{a}}\vk{B}_\alpha}
\left( \uth ^{-\frac{1}{\beta}}
+
 \oth ^{-\frac{1}{\beta}}
\right)
\Gamma\LT(\frac{1}{\beta}+1\RT)u^{\frac{2}{\alpha}-\frac{2}{\beta}} \prod_{i=1}^n\Psi\LT(c_i u\RT)\ooo.
\EQNY

\underline{Case ii) $\alpha=\beta$}: {\it Upper bound.} Bonferroni inequality implies
\BQNY
P_1(u)&\le&  \pk{\sup_{t\in D_u}\min_{1\le i\le n}Z_i^{\vn_+}(t)> u}
 \le  \pk{\sup_{t\in \Del_{-1}\cup\Del_0}\min_{1\le i\le n}Z_i^{\vn_+}(t)> u}+\Sigma_4(u)+\Sigma_5 (u),
\EQNY
where
\BQNY
\Sigma_4(u)= \sum_{ k=1}^{ N(u)} \pk{\sup_{t\in \Del_k}\min_{1\le i\le n}Z_i^{\vn_+}(t)> u},\ \ \Sigma_5 (u)=\sum_{ k=-N(u)}^{-2}\pk{\sup_{t\in \Del_k}\min_{1\le i\le n}Z_i^{\vn_+}(t)> u}.
\EQNY
It follows from \nekorr{lemPit} that
\BQNY
\pk{\sup_{t\in \Del_{-1}\cup\Del_0}\min_{1\le i\le n}Z_i^{\vn_+}(t)> u}=
\mathcal{H}_{\vk{c} \sqrt{\vk{a}^{\vn_+}}\vk{B}_\alpha, \vk{c}^2\vk{d}^{\vn_-}}[-T_1,T_1] \prod_{i=1}^n\Psi\LT(c_i u\RT)\ooo
\EQNY
as $u\to \IF$,
where $\vk{d}^{\vn_\pm}(t) = (d_1^{\vn_\pm}(t),\cdots,d_n^{\vn_\pm}(t))$ with
 $d_i^{\vn_\pm}(t)$ given as in \eqref{eq:d}.
Moreover, the same arguments as in the derivation of \eqref{eq:upp2}
yield
\BQNY
\Sigma_4(u)&\le& \mathcal{H}_{\vk{c} \sqrt{\vk{a}^{\vn_+}}\vk{B}_\alpha}(T_1)\sum_{k=1}^{N(u)}\prod_{i=1}^n\Psi(f_i^{\vn_-}(k,u))\ooo\\
&\le&
\mathcal{H}_{\vk{c} \sqrt{\vk{a}^{\vn_+}}\vk{B}_\alpha}(T_1 )  \prod_{i=1}^n\Psi\LT(c_i u\RT)
\sum_{k=1}^{\IF}\exp\LT(-\oth_{\vn-}  (kT_1)^\beta\RT)\ooo
\EQNY
and similarly
\BQNY
\Sigma_5(u)
 \le  \mathcal{H}_{\vk{c} \sqrt{\vk{a}^{\vn_+}}\vk{B}_\alpha}(T_1)  \prod_{i=1}^n\Psi\LT(c_i u\RT)
\sum_{k=1}^{\IF}\exp\LT(-   \uth_{\vn-}    (kT_1)^\beta\RT)\ooo.
\EQNY
{\it Lower bound.} Let $T_2$ be any positive constant. We have by \nekorr{lemPit} that
\BQNY
P_1(u)&\ge&    \pk{\sup_{t\in \LT[-T_2u^{-\frac{2}{\alpha}}, T_2u^{-\frac{2}{\alpha}}\RT]}\min_{1\le i\le n}Z_i^{\vn_-}(t)> u}\\
&=&\mathcal{H}_{\vk{c} \sqrt{\vk{a}^{\vn_-}}\vk{B}_\alpha, \vk{c}^2\vk{d}^{\vn_+}}[-T_2,T_2] \prod_{i=1}^n\Psi\LT(c_i u\RT)\ooo
\EQNY
as $u\to\IF$. From the upper and lower bounds obtained above we conclude that
\BQNY
\mathcal{H}_{\vk{c} \sqrt{\vk{a}^{\vn_-}}\vk{B}_\alpha, \vk{c}^2\vk{d}^{\vn_+}}[-T_2,T_2] &\le&
 \liminf_{u\to\IF}\frac{P_1(u)}{\prod_{i=1}^n\Psi\LT(c_i u\RT)}\\
&\le&
 \limsup_{u\to\IF}\frac{P_1(u)}{\prod_{i=1}^n\Psi\LT(c_i u\RT)}\\
&\le &\mathcal{H}_{\vk{c} \sqrt{\vk{a}^{\vn_+}}\vk{B}_\alpha, \vk{c}^2\vk{d}^{\vn_-}}[-T_1,T_1]
+
\mathcal{H}_{\vk{c} \sqrt{\vk{a}^{\vn_+}}\vk{B}_\alpha}(T_1)
\sum_{k=1}^{\IF}\exp\LT(-\uth_{\vn-} (kT_1)^\beta\RT)\\
&&+
\mathcal{H}_{\vk{c} \sqrt{\vk{a}^{\vn_+}}\vk{B}_\alpha}(T_1)
\sum_{k=1}^{\IF}\exp\LT(-\uth_{\vn+} (kT_1)^\beta\RT).
\EQNY
In the light of \eqref{subH}, using that
$
\uth>0,
\oth>0
$,
and letting
$\vn\to0, T_2\to\IF$
on the left-hand side of the last equation we
obtain that
$\mathcal{H}_{\vk{c} \sqrt{\vk{a}}\vk{B}_\alpha}^{\vk{c}^2\underline{\vk{b}},\vk{c}^2\overline{\vk{b}}}<\IF$. Similarly, letting $\vn\to0, T_1\to\IF$ on the right-hand side implies
$\mathcal{H}_{\vk{c} \sqrt{\vk{a}}\vk{B}_\alpha}^{\vk{c}^2\underline{\vk{b}},\vk{c}^2\overline{\vk{b}}}>0$.
Therefore, we conclude that
\BQNY
P_1(u) = \mathcal{H}_{\vk{c} \sqrt{\vk{a}}\vk{B}_\alpha}^{\vk{c}^2\underline{\vk{b}},\vk{c}^2\overline{\vk{b}}}
 \prod_{i=1}^n\Psi\LT(c_i u\RT)\ooo,\ \ u\to\IF.
\EQNY

\underline{Case iii) $\alpha>\beta$}:

{\it Upper bound.}
Since $\alpha>\beta$, we have that
\BQNY
P_1(u)
&\le&
\pk{\sup_{t\in D_u}\min_{1\le i\le n}Z_i^{\vn_+}(t)> u}
\le \pk{\sup_{t\in \Del_{-1}\cup\Del_0}\min_{1\le i\le n}Z_i^{\vn_+}(t)> u}\\
 &\le &\pk{\sup_{t\in \Del_{-1}\cup\Del_0}\min_{1\le i\le n}W_i^{\vn_+}(t)> u},
\EQNY
where
$$
W_i^{\vn_+}(t)=\frac{\sigma_{X_i}(t_0)}{1+d_i^{\vn_{-}}(t)}V_i(t),\ \ \ t\in\R,
$$
with $V_i(t)$ being a centered stationary Gaussian process with covariance
{$r_{V_i}(t)=\exp(-|t|^\beta), t\ge0$.}
Thus, by the proof of Case $ii)$
\BQNY
P_1(u)
 \le
\mathcal{H}_{\vk{c}\vk{B}_\beta, \vk{c}^2 \vk{d}^{\vn_{-}}}[-T_1,T_1]
\prod_{i=1}^n\Psi\LT(c_i u\RT)\ooo,\ \ u\to\IF.
\EQNY
{\it Lower bound.} It follows easily that
\BQNY
P_1(u) \ge   \pk{\min_{1\le i\le n}X_i(t_0)> u}
 =  \prod_{i=1}^n\Psi\LT(c_i u\RT).
\EQNY
Letting  $T_1\to0$  we have from the above upper and lower bounds that
\BQNY
P_1(u) =\prod_{i=1}^n\Psi\LT(c_i u\RT)\ooo,\ \ u\to\IF.
\EQNY

\underline{ Asymptotics of $P_2(u)$}:
In order to complete the proof we show that \eqref{eq:P2} is valid.  To this end, we shall derive an adequate upper bound for $P_2(u)$  by utilizing the generalized Borell-TIS and Piterbarg inequalities given in \nelem{GBorell} and \nelem{GPiter}, respectively.

By Assumption I and Assumption III we can choose some small $\vn>0$ such that
\BQNY
\inf_{t\in[0,T] \LE{\setminus}[t_0-\vn,t_0+\vn]}g(t)\ > \ g(t_0).
\EQNY
Clearly
\BQNY
 P_{2}(u)\le \pk{\sup_{t\in[0,T] \LE{\setminus} [t_0-\vn,t_0+\vn]}\min_{1\le i\le n}X_i(t)> u } +\pk{\sup_{t\in[t_0-\vn,t_0+\vn]\LE{\setminus}(t_0+D_u)}\min_{1\le i\le n}X_i(t)> u }.
\EQNY
Further, by the generalized Borell-TIS  \ccc{inequality} in \nelem{GBorell} for all sufficiently large $u$ we have
\BQN \label{eq:asy1}
\pk{\sup_{t\in[0,T] \LE{\setminus} [t_0-\vn,t_0+\vn]}\min_{1\le i\le n}X_i(t)> u }\le \exp\LT(-\frac{(u-\mu)^2}{2} \inf_{t\in[0,T]\LE{\setminus}[t_0-\vn,t_0+\vn]}g(t)\RT)
\EQN
for some $\mu>0$. Moreover, in view of the generalized Piterbarg inequality in \nelem{GPiter} we obtain for all sufficiently large $u$
\BQNY
\pk{\sup_{t\in[t_0-\vn,t_0+\vn]\LE{\setminus}(t_0+D_u)}\min_{1\le i\le n}X_i(t)> u }\le Q_1 u^{\frac{2}{\min(\gamma,\alpha)}-1}\exp\LT(-\frac{u^2}{2} \inf_{t\in[t_0-\vn,t_0+\vn]\LE{\setminus}(t_0+D_u)}g(t)\RT)
\EQNY
for some $Q_1>0$. In addition (recall that
$\uth>0$,
$\oth>0$),
we have that
$$
g(t_0+t)\ge g(t_0)+Q_2 (\delta(u))^\beta
$$
holds for all $t\in[ -\vn, \vn] \LE{\setminus}D_u$. Therefore, for all $u$ large
\BQN\label{eq:asy2}
\pk{\sup_{t\in[t_0-\vn,t_0+\vn]\LE{\setminus}D_u}\min_{1\le i\le n}X_i(t)> u }\le Q_1 u^{\frac{2}{\min(\gamma,\alpha)}-1}\exp\LT(-\frac{u^2}{2} g(t_0)\RT)
\exp\LT(-\frac{Q_2}{2}(\ln u)^2\RT).
\EQN
Consequently, the claim in \eqref{eq:P2} follows immediately from \eqref{eq:asy1}, \eqref{eq:asy2} and the asymptotics of $P_1(u)$ obtained above. Thus the proof is complete. \QED 

\bigskip

\textbf{Acknowledgments.}  We are thankful to  the associate editor \ccc{and the referees} for several suggestions which improved our manuscript. Partial support from an SNSF grant \tbb{and RARE -318984 (an FP7 Marie Curie IRSES Fellowship) is kindly acknowledged.}
The first author also acknowledges  partial support by NCN Grant No 2013/09/B/ST1/01778 (2014-2016).

\bibliographystyle{plain}
\bibliography{sYZ}
\end{document}